\documentclass[11pt]{article}
\title{Stochastic grid bundling method for backward stochastic differential equations}

\usepackage{fullpage}
\usepackage{amsfonts,amssymb,amsmath,amsthm}
\usepackage{array}
\usepackage{graphicx}
\usepackage{caption}
\usepackage{subcaption}
\usepackage{color}
\usepackage{multirow}
\usepackage{longtable}
\usepackage{authblk}

\usepackage{tikz}
\usetikzlibrary{shapes.geometric, arrows}

\author[1]{Ki Wai Chau}
\author[1,2]{Cornelis W. Oosterlee}
\affil[1]{Research Group of Scientific Computing, Centrum Wiskunde \& Informatica}
\affil[2]{Department of Applied Mathematics, Delft University of Technology}

\allowdisplaybreaks[1]
\numberwithin{equation}{section}

\newtheorem{theorem}{Theorem}[section]

\newtheorem{proposition}[theorem]{Proposition}
\theoremstyle{definition}\newtheorem{definition}[theorem]{Definition}
\theoremstyle{definition}\newtheorem{assumption}[theorem]{Assumption}
\theoremstyle{remark}\newtheorem{remark}[theorem]{Remark}

\newcommand{\expectation}[3]{\mathbb{E}^{#3}_{t_{#2}}\left[#1\right]}
\newcommand{\transpose}[1]{#1^T}
\newcommand{\approximant}[3]{#1^{#3}_{#2}}
\newcommand{\regressionparameter}[2]{#1_{#2}(b)}
\newcommand{\leastsquare}[3]
{
	\arg\min_{#1 \in \mathbb{R}^Q} \frac{\sum^{M}_{m = 1}
	(p(X^{\pi,m}_{t_{#2+1}})#1-#3)^2{\bf 1}_{\mathcal{B}_{t_k}(b)}(X^{\pi, m}_{t_k})}{\sum^M_{m=1}{\bf 1}_{\mathcal{B}_{t_k(b)}}(X^{\pi. m}_{t_k})}\
}
\newcommand{\itemassumption}[1]{$({\bf A}_{#1})$}
\newcommand{\path}[2]{X^{\pi, #2}_{t_{#1}}}

\begin{document}

\maketitle

\begin{abstract}
In this work, we apply the Stochastic Grid Bundling Method (SGBM) to numerically solve backward stochastic differential equations (BSDEs). 
The SGBM algorithm is based on conditional expectations approximation by means of bundling of Monte Carlo sample paths and a local regress-later regression within each bundle.
The basic algorithm for solving the backward stochastic differential equations will be introduced and an upper error bound is established for the local regression.
A full error analysis is also conducted for the explicit version of our algorithm and numerical experiments are performed to demonstrate various properties of our algorithm.
\end{abstract}

\section{Introduction}\label{section_introduction}

The Stochastic Grid Bundling Method (SGBM) is a Monte Carlo based algorithm designed to solve backward dynamic programming problems, with applications in pricing Bermudan options in \cite{jain_oosterlee_2015} and \cite{cong_oosterlee_2015}.
This algorithm has been further extended computationally by the incorporation of GPU acceleration in \cite{leitao_oosterlee_2015} and generalized to the computation of Credit Valuation Adjustment and Potential Future Exposure in \cite{degraaf_etc_2014}.
In this work, we will extend its applicability to the approximation of Backward Stochastic Differential Equations (BSDEs). 
We shall also study the  errors in the SGBM algorithm.

The SGBM algorithm is based on the so-called {\em regress-later technique} and on an adaptive local basis approach.
In usual Monte Carlo regression methods for backward-in-time problems, the values of the target function at the end of a time interval are regressed on certain dependent variables that are measured at the beginning of the time interval (which is called the regress-now approach). This creates a statistical error.
Instead, the dependent variable is projected onto a set of basis functions at the end of the interval in a regress-later method, and a conditional expectation across the interval is then computed for each basis function.
This difference removes the statistical error in the regression step.
Regress-later schemes have been further discussed in \cite{Glasserman2004}.

With an adaptive local basis approach, the whole simulation is partitioned into non-overlapping subsets and we perform least-squares regressions separately within these subsets, possibly with a different basis for each subset.
The exact partition depends on the simulated samples themselves and its purpose is to gather samples that share similar "characteristics" such that the local regression is more accurate than the global one.
For further application of localization in numerical schemes, the reader may check out \cite{Bouchard2012}.
Since each partition is non-overlapping, SGBM is easy to scale up in dimensionality and can facilitate parallel computing.
We would like to test the SGBM algorithm in a new problem setting such that we can take advantage of its nice properties and also get a better understanding of the underlying principles.

The problem that we are interested in is the numerical approximation of BSDEs. 
These equations form a popular subject of research in quantitative finance ever since their introduction in \cite{pardoux_peng_1992} and related works.
The connection between BSDEs and partial differential equations (PDEs) also provides the opportunity of solving PDEs (in high dimensions) with stochastic methods.
However, the computational difficulties of solving BSDEs prevent them from being widely used in practice.
Therefore, efficient algorithms for the approximation of high-dimensional BSDEs are of great interest.
In fact, there are numerous works just focusing on the application of Monte Carlo methods to BSDEs, including \cite{BOUCHARD2004175, CRISAN20101133, doi:10.1137/16M106371X, lemor2006, Bender2012} and some of these could be integrated with our proposed scheme for further development.
For example, in \cite{ding2017aregression}, the authors proposed a regression basis based on a Fourier-cosine expansion in a least-squares scheme for BSDEs, which can possibly be used as a basis function in our SGBM algorithm.
As far as we know, there is no study of a combined approach based on the regress-later scheme and a localization Monte Carlo technique for these equations, which is the goal of this work.

In this article, we consider the application of SGBM to decoupled Forward Backward Stochastic Differential Equations (FBSDEs) of the form 
\begin{equation}
\label{equation_fbsde}
\left\{
\begin{array}{l}
dX_t = \mu(t,X_t)dt + \sigma(t, X_t)dW_t, \; X_0 = x_0,\\
dY_t = -f(t, X_t, Y_t, Z_t)dt + Z_tdW_t, \; Y_T = \Phi(X_T),
\end{array}
\right.
\end{equation}
defined on $0 \leq t \leq T$.
The function $f:[0,T]\times \mathbb{R}^q \times \mathbb{R} \times \mathbb{R}^d$ is called the driver function of the backward process and the process  $W_t = (W_{1,t}, \ldots, W_{d,t})^\top$ is a d-dimensional standard Brownian motion.
Note that the usual setting of complete probability space $(\Omega, \mathcal{F},\mathbb{F},\mathbb{P})$ with $\mathbb{F}:=(\mathcal{F}_t)_{0\leq t\leq T}$ being a filtration satisfying the usual conditions of being right-complete and $\mathbb{P}$-complete applies throughout the article.
Given that a solution exists for the forward equation, a pair of adapted processes $(Y_t,Z_t)$ is said to be the solution of the FBSDE \eqref{equation_fbsde}, if $Y_t$ is a continuous real-valued adapted process, $Z_t$ is a real-valued predictable process such that $\int^T_0|Z_t|^2dt < \infty$ almost surely in $\mathbb{P}$ and the pair satisfies Equation \eqref{equation_fbsde}.

One key difficulty in solving a BSDE is that the pair $(Y_t, Z_t)$ must be adapted to the underlying filtration.
The terminal condition $Y_T$ is given by $\Phi(X_T)$, where $\Phi:\mathbb{R}^q\rightarrow\mathbb{R}$ is a deterministic function. 
Therefore, $Y_T$ is adapted to the filtration $\mathcal{F}_T$ and a naive Euler discretization on the backward equation fails to produce an adapted solution, for further discussion on this, the reader may refer to the introduction in \cite{BOUCHARD2004175}.
In this work, we aim to construct an approximate solution by the theta-scheme from \cite{zhao_wang_peng_2009} and applying the SGBM algorithm.

To ensure the existence and uniqueness of the solution to the forward equation, further regularity conditions are assumed here.
The functions $\mu: [0,T] \times \mathbb{R}^q\rightarrow \mathbb{R}^q$ and $\sigma: [0,T] \times \mathbb{R}^q \rightarrow \mathbb{R}^{q \times d}$ refer to the drift and the diffusion coefficients of the forward stochastic process, and $x_0$ is the initial condition for $X$.  
It is assumed that both $\mu(t,x)$ and $\sigma(t,x)$ are measurable functions that are uniformly Lipschitz in $x$ and such condition holds uniformly in $t$. 
The forward process also satisfies the Markov property, namely $\mathbb{E}[X_\tau|\mathcal{F}_t] = \mathbb{E}[X_\tau|X_t]$ for $\tau \geq t$, where $\mathbb{E}[\cdot]$ denotes expectation with respect to probability measure $\mathbb{P}$.

The rest of the article is organized as follows.
We start in Section \ref{section_algorithm} with the introduction of the SGBM algorithm, along with the necessary time discretization scheme and assumptions. 
Section \ref{section_refined_regression} will present an error analysis of a simplified case of SGBM.
The proof in this section forms the foundation for the error bound in any algorithm applying SGBM.
Later, in Section \ref{section_explicit}, we derive the full error bound for a specific choice of discretization scheme as an example.
The article finishes with numerical experiments and a conclusion.

To close off this section, here is some further notation that is used in this article.
\begin{itemize}
	\item For any vector $x$, $|x|$ denotes its Euclidean norm and $x_r$ denotes its $r$-th component.
	\item Similarly, $X_{r,t}$ denotes the $r$-th component for any random process $X_t$.
	\item The gradient $\nabla g$ is defined as $\left(\frac{\partial g}{\partial x_1}, \ldots , \frac{\partial g}{\partial x_q}\right)$ for any differentiable function $g:\mathbb{R}^q\rightarrow \mathbb{R}$.
	\item The notations $\mathbb{E}_{t}[\cdot]$ and $\mathbb{E}^x_t[\cdot]$ are the simplified notations for $\mathbb{E}[\cdot|\mathcal{F}_t]$ and $\mathbb{E}[\cdot|X_t = x]$
	\item For any set $\mathcal{S}$, the function ${\bf 1}_S$ is the indicator function which takes value $1$ when the input is within set $\mathcal{S}$ and $0$ otherwise.  
	\item For any function space $H$ containing functions $\phi: \mathbb{R}^q \rightarrow \mathbb{R}$, $H^+$ is defined as the set $\{\{(x, y)\in \mathbb{R}^q \times \mathbb{R} : \phi(x) \geq y\} : \phi \in H\}$.
	\item For any function $\phi$ and compact set $\mathcal{A}$, the control constant $C_{\phi, \mathcal{A}}$ is defined as an extended real number $\sup_{x \in \mathcal{A}}|\phi(x)|$.
\end{itemize}

\section{Assumptions and Algorithm} \label{section_algorithm}

In this section, we shall introduce the SGBM algorithm and its application to the approximation of BSDEs.
To begin, we need to discretize the BSDEs.

\subsection{Discretization Scheme}\label{section_discretization scheme}

We denote a time grid $\pi=\{0 = t_0 < \ldots < t_N = T\}$ on the interval $[0,T]$ and let $\Delta_k :=t_{k+1}-t_k$, $\Delta W_{l,k} := W_{l,t_{k+1}}-W_{l,t_k}$, and $\Delta W_k := (\Delta W_{1,k},\ldots, \Delta W_{d,k})^\top$ be the time-step, the Brownian motion increment along the $l$-th dimension and the Brownian motion increment, respectively, for $k \in \{0, \ldots, N-1\}$.

For the forward process $X_t$, we shall apply a Markovian approximation $X^\pi_{t_k}, t_k\in \pi$.
The most common choice is the Euler-Maruyama scheme, which will be explained in Section \ref{section_numerical}.
However, our algorithm can work with any simulation method where the conditional expectations over one time step are known for some specific functions.

The backward in time discretizations $(Y^{\pi}, Z^{\pi})$ form a special case of the theta-scheme from \cite{zhao_li_zhang_2012} and \cite{ruijter_oosterlee_2015} by selecting $(\theta_1, \theta_2) = (0,1)$. 
We use the following explicit discretization:
\begin{align*}
& y_N(x) = \Phi(x),\\
& z_k(x) = \frac{1}{\Delta_k}\expectation{y_{k+1}(X^\pi_{t_{k+1}})\Delta W_k }{k}{x} , \;
k = N-1, \ldots, 0,\\
& y_k(x) = \expectation{y_{k+1}(X^\pi_{t_{k+1}})}{k}{x} + \Delta_k \expectation{f_{k+1}(y_{k+1}(X^\pi_{t_{k+1}}),z_{k+1}(X^\pi_{t_{k+1}}))}{k}{x},\; k= N-1,\ldots,0,
\end{align*}
where $f_k(y, z) := f(t_k, X^\pi_{t_k}, y, z)$.

\subsection{Standing Assumptions}\label{session_assumptions}

To ensure the existence and uniqueness of the solution of the continuous BSDEs, some basic assumptions are required.
Moreover, these assumptions will affect the algorithm designed regarding the admissible choice of $\pi$ and the error bound of the scheme.
In this work, we assume the global Lipschitz condition as stated in Assumption \ref{assumption:globally_lipschitz}.
Note that this assumption will affect the derivation and the result of the error bound for the complete algorithm.
Assumption \ref{assumption:globally_lipschitz} is in force here as it is the most common assumption in the BSDE literature. 
Alternative assumptions can be found, for instance, in \cite{gobet_turkedjiev_2016}.

\begin{assumption}[Globally Lipschitz driver] \label{assumption:globally_lipschitz} \quad\\
	\begin{enumerate} 
		\item[\itemassumption{\xi}]
		\begin{enumerate}
			\item[i.)]	$\Phi$ is a measurable function.
			\item[ii.)] The control constant $C_{\Phi, \mathcal{A}} < \infty$ for any given compact set $\mathcal{A}$.
		\end{enumerate}	 
		\item[\itemassumption{F}]
		\begin{enumerate}
			\item[i.)]
			$(t, x, y, z) \mapsto f(t, x, y, z)$ is $\mathcal{B}(\mathcal{\mathbb{R}}) \otimes \mathcal{B}(\mathbb{R}^q) \otimes \mathcal{B}(\mathbb{R}) \otimes \mathcal{B}(\mathbb{R}^d)$-measurable.
			\item[ii.)] For every $k \leq N$, $f_k(y,z)$ as defined in the Subsection \ref{section_discretization scheme} is $\mathcal{F}_{t_k}\otimes \mathcal{B}(\mathbb{R}) \otimes \mathcal{B}(\mathbb{R}^d)$-measurable and there exists an $L_f \in [0, +\infty )$ such that 
			\begin{equation*}
			|f_k(y,z)-f_k(y', z')| \leq L_f(|y-y'|+|z-z'|),\qquad\forall k \in \{0, \ldots, N\},
			\end{equation*}
			for any $(y, y', z, z') \in \mathbb{R} \times \mathbb{R} \times \mathbb{R}^d \times \mathbb{R}^d$.
			\item[iii.)] There exists a $C_f \in [0, \infty)$ such that 
			\begin{equation*}
			|f_k(0,0)| \leq C_f,\qquad \forall k \in \{0, \ldots, N\}.
			\end{equation*}
			\item[iv.)] The time discretization is such that 
			\begin{equation*}
			\limsup_{N\rightarrow\infty} R_\pi < +\infty, \quad \mbox{where } R_\pi = \sup_{0\leq k\leq N-2}\frac{\Delta_k}{\Delta_{k+1}}.
			\end{equation*}
		\end{enumerate}
	\end{enumerate}
\end{assumption}

Again, the assumption here is for the consistency of our derivation and does not imply that our algorithm can only be applied when these assumptions are satisfied.

\subsection{Stochastic Grid Bundling Method}

We now introduce SGBM.
Due to the Markovian setting of $(X^\pi_{t_k},\mathcal{F}_{t_k})_{t_k\in\pi}$, there exist functions $y_k(x)$ and $z_k(x)$ such that 
\begin{equation*}
Y^\pi_{t_k}=y_k(X^\pi_{t_k}),\; Z^\pi_{t_k} = z_k(X^\pi_{t_k}).
\end{equation*}
Our method is based on estimating these functions $(y_k(x),z_k(x))$ recursively backward in time by a local least-squares regression technique onto a finite function space with basis functions $(p_l)_{0\leq l \leq Q}$.

As a Monte Carlo based algorithm, our program starts with the simulation of $M$ independent samples of $(X^\pi_{t_k})_{0 \leq k \leq N}$, denoted by $(X^{\pi,m}_{t_k})_{1 \leq m \leq M, 0 \leq k \leq N}$.
Note that in this basic algorithm, the simulation is only performed once. 
This scheme is therefore a non-nested Monte Carlo scheme.

The next step is the backward recursion. 
Denote by $y^{R}_k$ the SGBM approximation of the function $y_k$.
The function $\approximant{z}{k}{R}$ similarly means the approximation of $\approximant{z}{k}{}$.

At initialization, we assign the terminal values to each path for our approximations, i.e., 
\begin{equation*}
\approximant{y}{N}{R}(X^{\pi, m}_{t_N})=\Phi(X^{\pi, m}_{t_N}),\quad 
m = 1,\ldots, M,
\end{equation*}
The following steps are performed recursively, backward in time, at $t_k$, $k=N-1,\ldots, 0$.
First, we bundle all paths into $\mathcal{B}_{t_k}(1),\ldots,\mathcal{B}_{t_k}(B)$ non-overlapping partitions based on the result of $(X^{\pi, m}_{t_k})$.
Note that our design allows the application of various clustering techniques within the SGBM algorithm. 
A previous study in \cite{leitao_oosterlee_2015} compares the k-means clustering with an equal partitioning, and shows that they are similar in accuracy.
However, it remains an interesting problem which clustering technique would provide the optimal result.
We use the equal partition technique, which will be specified in Section \ref{section_refined_regression}, for the error analysis and the numerical experiment.

Next, we perform the regress-later approximation separately within each bundle.
The regress-later technique we are using combines the least-squares regression with the (analytical) expectations of the basis functions to calculate the necessary expectations.

Generally speaking, for $M$ Monte Carlo paths, a standard regress-now algorithm for a dynamic programming problem finds a function $\iota$ within the space spanned by the regression basis such that it minimizes the value  $\frac{1}{M}\sum^M_{i=1}(g(X^i_{t+\delta}) - \iota(X^i_{t}))^2$ and approximates the expectation $\mathbb{E}_t[g(X_{t+\delta})]$ by $\mathbb{E}_t[\iota(X_t)] = \iota(X_t)$.
As a projection from a function of $X_{t+\delta}$ to a function of $X_t$ is performed then, it would introduce a statistical bias to the approximation.

Instead, the regress-later technique we employ picks out a function $\kappa$ such that it minimizes  $\frac{1}{M}\sum^M_{i=1}(g(X^i_{t+\delta}) - \kappa(X^i_{t+\delta}))^2$ and approximates the expectation $\mathbb{E}_t[g(X_{t+\delta})]$ by $\mathbb{E}_t[\kappa(X_{t+\delta})]$.
By using functions on the same variable in the regression basis, we can remove the statistical bias in the regression.
However, the expectation of all basis functions must preferably be known in order to apply the regress-later technique efficiently.

In the context of our algorithm, we define the bundle-wise regression parameters $\regressionparameter{\alpha}{k+1}$, $\regressionparameter{\beta}{k+1}$, $\regressionparameter{\gamma}{k+1}$ as 
\begin{align*}
& \regressionparameter{\alpha}{k+1} = \leastsquare{\alpha}{k}{\approximant{y}{k+1}{R}(\path{k+1}{m})},\\
& \regressionparameter{\beta}{i, k+1} = 
\leastsquare{\beta}{k}{\approximant{z}{i, k+1}{R}(\path{k+1}{m})},\\
& \regressionparameter{\gamma}{k+1} = 
\leastsquare{\gamma}{k}{f_{k+1}(\approximant{y}{k+1}{R }(\path{k+1}{m}), \approximant{z}{k+1}{R}(\path{k+1}{m}))}.
\end{align*}
The approximate functions within the bundle at time $k$ are defined by the above parameters and the expectations $\mathbb{E}^x_{t_k}[p(X^\pi_{t_{k+1}})]$ and $\mathbb{E}^x_{t_k}\left[p(X^\pi_{t_{k+1}})\frac{\Delta W_{r, k}}{\Delta_k}\right]$:
\begin{align*}
& \approximant{z}{r, k}{R}(b, x) =  \expectation{
	\frac{\Delta W_{r,k}}{\Delta_k}p(X^\pi_{t_{k+1}})
}{k}{x} \regressionparameter{\alpha}{k+1},\quad
r= 1,\ldots,d;\\
& \approximant{y}{k}{R}(b, x) =
\expectation{p(X^\pi_{t_{k+1}})}{k}{x}
(\regressionparameter{\alpha}{k+1}+ \Delta_k(1 - \theta_1)\regressionparameter{\gamma}{k+1}),
\quad i = 1, \ldots, I.
\end{align*}

As the expectations related to the basis functions are the foundation of any regress-later scheme, we assume that the following assumptions are satisfied. 
\begin{assumption}
	The regression basis $\{p_1, \ldots, p_Q\}$ is assumed to satisfy the following assumptions.
	\begin{enumerate} 
		\item[\itemassumption{p}] 
		\begin{enumerate}
			\item[i.)]
			$\mathbb{E}^x_{t_k}[p_l(X^\pi_{t_{k+1}})]$ and  $\mathbb{E}^x_{t_k}\left[p_l(X^\pi_{t_{k+1}})\frac{\Delta W_{r, k}}{\Delta_k}\right]$ are known, either analytically or empirically, for all $k = 0, \ldots N-1$, $l = 1, \ldots, Q$ and $r = 1, \ldots, d$.
			\item[ii.)] For any given compact set $\mathcal{A}$ in $\mathbb{R}^q$, the constant $C_{p, \mathcal{A}} := \max_{l=1, \ldots, Q}C_{p_l, \mathcal{A}}$. Moreover, there exists a constant  $C_{M,\mathcal{A}}$ such that 
			\begin{equation*}
				\sum^Q_{l=1} \left|\mathbb{E}^x_{t_k}[p_l(X^\pi_{t_{k+1}})]\right| \leq C_{M, \mathcal{A}},\qquad \forall x \in \mathcal{A}, \mbox{ and } k = 0, \ldots, N-1;
			\end{equation*} 
			and 
			\begin{equation*}
				\sum^Q_{l=1} \left|\mathbb{E}^x_{t_k}[p_l(X^\pi_{t_{k+1}})\frac{\Delta W_{r, k}}{\Delta_k}]\right| \leq C_{M,\mathcal{A}},\qquad \forall x \in \mathcal{A}, \mbox{ and } k = 0, \ldots, N-1.
			\end{equation*}
		\end{enumerate}
	\end{enumerate}
\end{assumption}

Next, to ensure the stability of our algorithm, $|\alpha_k(b)|$, $|\beta_{r, k}(b)|$ and $|\gamma_k(b)|$ must be bounded above for all $k, b, r$.
In practice, this means that an error notion should be given by the program when the Euclidean norm of any regression coefficient vector is greater than a predetermined constant $L$.
Further details on this requirement will be described in Section \ref{section_refined_regression}.

Finally, to simplify notation, we define the notations below for the regression result across the bundles.
\begin{align*}
\tilde{y}_{k+1}^{R, I}(x_1, x_2) & := \sum_{b = 1}^B{\bf 1}_{\mathcal{B}_{t_k}(b)}(x_1)p(x_2)\alpha_{k+1}(b),\\
\tilde{z}_{r, k+1}^{R} (x_1, x_2) & := \sum_{b = 1}^B{\bf 1}_{\mathcal{B}_{t_k}(b)}(x_1)p(x_2)\beta_{r, k+1}(b),\\
\tilde{f}_{k+1}^{R}(x_1, x_2) & := \sum_{b = 1}^B{\bf 1}_{\mathcal{B}_{t_k}(b)}(x_1) p(x_2)\gamma_{k+1}(b).
\end{align*}

\section{Refined Regression}\label{section_refined_regression}
In this section, we derive a proof of an error bound for our regress-later strategy.
In order to ensure the stability of our algorithm, we have introduced a sample selection step into the algorithm and modified the classical proof for nonparametric regression from \cite{gyorfi_kohler_krzyzak_walk_2002}, which was used in \cite{gobet_turkedjiev_2016}, for the derivation of the error bound to SGBM.

In order to simplify expressions, different notations are used in this section.
We consider a random vector $(X,Y)$, where $X$ and $Y$ are both $\mathbb{R}^q$, following the probability measure $\nu$.
A cloud of simulation paths can be generated by independently simulating $M$ copies, $\{(X^m, Y^m) : m = 1, \ldots, M\}$, defined on a probability space $(\hat{\Omega}, \hat{\mathcal{F}}, \hat{\mathbb{P}})$.
In our content, the pair $(X, Y)$ represents the independent and dependent variables under consideration and $(X^m, Y^m)$ are the simulated samples for $(X, Y)$.

Denote by $\mathbb{B}$ a specific partition with $\mathbb{B} :=\{\mathcal{B}(1), \ldots, \mathcal{B}(B)\}$ and $\bigcup_{b=1}^B\mathcal{B}(b) = \mathbb{R}^d$.
The partition which is used in the regression estimates is based on the simulation data $X^m$ in our setting and  to which bundle a sample belongs solely depends on $X^m$.

The main goal of SGBM is finding an effective and accurate way to approximate the expectation $\mathbb{E}\left[\left.v(Y)\right|X\right]$ in a recurrence setting for some deterministic function $v: \mathbb{R}^q \rightarrow \mathbb{R}$, and we begin with establishing an estimate $\tilde{v}: \mathbb{R}^q \times \mathbb{R}^q$ for $v$.
Note that although $v$ solely depends on the dependent variables, the estimate $\tilde{v}$ depends on both the independent and dependent variables in preparation for further calculation.

For a given partition and samples, one way to define the estimate $\tilde{v}$ is
\begin{align}
	\tilde{v}(x,y) 
	:= & 
	\sum^B_{b=1}{\bf 1}_{\mathcal{B}(b)}(x)\tilde{v}_b(y)
	=
	\sum^B_{b=1}{\bf 1}_{\mathcal{B}(b)}(x)\sum^{Q}_{k=1}\alpha_k(b) p_k(y),
	\label{equation_approximant}
\end{align}
where 
$$ \tilde{v}_b := \arg\min_{\phi \in H} \left\{ \frac{\sum^M_{m=1}{\bf 1}_{\mathcal{B}(b)}(X^m)|v(Y^m)-\phi(Y^m)|^2}{\sum^M_{m=1}{\bf 1}_{\mathcal{B}(b)}(X^m)}\right\}.$$
	
\begin{remark}
	It is possible that under some particular clustering scheme for SGBM, there would be empty bundles in the resulting partition. 
		
	In practice, one could simply ignore these empty bundles in the algorithm.		
	As there are no samples in these bundles, approximations within these bundles are not needed for the next time step. 
	One point to note is that since least-squares regression requires a sufficient number of samples to be accurate, adopting a bundling scheme that would produce a small number of samples in any bundle may not be a good idea. 
		
	When generalizing the theoretical proof below to bundling methods other than equal partition, one has to take this into account and define the measurable partition $\mathbb{B}$ in such a way that it is consistent with the practical bundling scheme while it also merges all empty bundle to non-empty ones.
		
	Further discussion on bundles with few samples under the equal partition scheme is placed in Remark \ref{remark_empty_bundle_equal}.
\end{remark}

Note that functions $\tilde{v}_b:\mathbb{R}^q\rightarrow\mathbb{R}$ are stochastic with respect to the simulation samples $(X^m, Y^m)$.
The linear vector space $H$ is spanned by continuous functions $\{p_1,\ldots, p_Q\}$, with $p_l: \mathbb{R}^q \rightarrow \mathbb{R}, \;\forall l = 1, \ldots, Q$.
Thus, the second equality in Equation \eqref{equation_approximant} just follows from the definition of $H$ and typical least-squares regression. 
In fact, if we denote the total number of samples in a given bundle by $\# \mathcal{B}(b)$ and let $\{(X^{b,1}, Y^{b,1}), \ldots, (X^{b,{\#\mathcal{B}(b)}}, Y^{b,{\#\mathcal{B}(b)}})\}$ be the samples in this bundle, the coefficients $\alpha(b)$ satisfy 
\begin{equation}
\mathcal{I}^\top\mathcal{I}\alpha(b) = \mathcal{I}^\top v(Y^b), \label{equation_coefficients}
\end{equation}
with $$\mathcal{I} = (p_j(Y^{b,i}))_{1\leq i \leq \#\mathcal{B}(b), 1 \leq j \leq Q} \mbox{ and } v(Y^b) = (v(Y^{b,1}), \ldots, v(Y^{b, \#\mathcal{B}(b)}))^\top.$$ 
According to \cite{gyorfi_kohler_krzyzak_walk_2002}, system \eqref{equation_coefficients} is always solvable, in the next section, we also provide a heuristic argument for its invertibility.
Again, the coefficients $\alpha$ in each bundle can be seen as random variables with respect to $(X^m, Y^m)_{m = 1, \ldots, M}$.

Reversely, we may select the simulation cloud based on the regression coefficients.  
Let the set $S$ be the set containing all possible collections of $(x_m, y_m)_{1 \leq m \leq M} \in (\mathbb{R}^q\times \mathbb{R}^q)^M$ such that $|\alpha(b)|^2 \leq L$ for all $b$ given that $(X^m, Y^m) = (x_m, y_m)_{1 \leq m \leq M}$.
We modify the probability of the simulation cloud by only accepting those results that are in $S$. 
We denote the modified expectation by $\hat{\mathbb{E}}_S$ and it is related to the original expectation by $\hat{\mathbb{E}}_S[{\bf 1}_A] = \frac{\hat{\mathbb{E}}[{\bf 1}_A{\bf 1}_S]}{\hat{\mathbb{E}}[{\bf 1}_S]}$.
\footnote{The situation of $\frac{0}{0}$ should be understood as $0$ and $\frac{K}{0}$ as $\infty$ in the rest of this article.}

\begin{remark}
	In a regress-now scheme, especially in a recursion scheme, the resulting approximation is truncated such that its value is within a bounded interval $[M_1, M_2]$.
	The truncation guarantees the convergence and the stability of the scheme. 
	However, truncation is not feasible in our regress-later scheme as we have to keep the full function for further operation. 
	Therefore, we must instead control the output by limiting the admissible samples. 	
\end{remark}

\begin{remark}
	The introduction of bundling here essentially serves two purposes.
	First of all, clustering data may act as a localization of function $v$, thus a more accurate approximation for $v$ can be achieved with lower order function basis. 
	This is especially beneficial for the high-dimensional case as basis functions in higher dimension are generally complicated and hard to calculate. 
	We need a method to increase accuracy without adding more basis functions.
	Secondly, by partitioning data into non-overlapping bundles, we can facilitate the application of parallel computing, which is important when we are in a high-dimensional situation.
	However, while the above benefit depends on the particular choice of basis, the analysis we do in this section is applicable for a more general setting. 
	So, we would not emphasis these points further in this section.
\end{remark}

Using $\nu$ to denote the probability measure induced by the random variable $(X, Y)$ and \\ $(X^m, Y^m)_{m = 1, \ldots, M}$ are independent and identical copies following the same law under a different probability space, the following random norms (depending on the simulation cloud $(X^m, Y^m)$) are used to quantify the error of approximation.
\begin{definition}
	Let $\varphi: \hat{\Omega} \times \mathbb{R}^q \times \mathbb{R}^q \rightarrow \mathbb{R}$ be measurable.
	For any set $\mathcal{B}\subset \mathbb{R}^q$, we define the following random norms
	\begin{equation*}
		||\varphi||^2_{\mathcal{B},\infty} := \frac{\int_\mathcal{B}\int|\varphi(x, y)|^2\nu(dx, dy)}{\int_\mathcal{B}\int\nu(dx, dy)}; \qquad ||\varphi||^2_{\mathcal{B}, \#} := \frac{\sum_{m=1}^M{\bf 1}_\mathcal{B}(X^m)|\varphi(X^m, Y^m)|^2}{\sum_{m=1}^M{\bf 1}_\mathcal{B}(X^m)}.
	\end{equation*}
\end{definition}
We derive the following theorem for the estimation of the error.
Since we only accept a simulation result that satisfies event $S$, we should only consider the average error among all these accepted events.
\begin{theorem} \label{theorem_regression}
	Assume that we perform an equal partition at the bundling step, namely, we order all samples according to some specific measurable sorting function on $X$, and separate them into almost-equal size bundles by the ordering.
	Further, assume a compact set $\mathcal{A} \subset \mathbb{R}^q$ to be given such that $C_{v, \mathcal{A}} \leq \infty$ and $\int v^2(y)\nu(dx, dy) \leq \infty$, namely, the function $v$ is within the $L^2$ space with respect to the given probability measure. 
	Then, for any real function $v$, we have
	\begin{align*}
		& \hat{\mathbb{E}}_S\left[\iint |v(y) - \tilde{v}(x,y)|^2 \nu(dx,dy)\right]\\
		\leq & 
		\frac{\vartheta(L')}{\hat{\mathbb{E}}[{\bf 1}_S]}\hat{\mathbb{E}}
		\left[\sum^B_{b = 1}\int_{\mathcal{B}(b)}\int\nu(dx, dy)\frac{(\log(\sum^M_{m=1}{\bf 1}_{\mathcal{B}(b)}(X^m) - 1)+1)(Q+1)}{\sum^M_{m=1}{\bf 1}_{\mathcal{B}(b)}(X^m) - 1}\right]\\
		& 
		+ \hat{\mathbb{E}}_S\left[\sum^{B-1}_{b =1}\int_{\mathcal{B}(b)}\int\nu(dx, dy)\frac{24 L'}{(\sum_{m=1}^M{\bf 1}_{\mathcal{B}(b)}(X^m))}\right]\\
		& 
		+ \frac{12}{\hat{\mathbb{E}}[{\bf 1}_S]}
		\hat{\mathbb{E}}\left[
		\sum_{\mathcal{B}\in \mathbb{B}}\int_\mathcal{B}\int\nu(dx,dy)
		(\inf_{\phi \in H} \sup_{x\in \mathcal{B}}\mathbb{E}\left[ |v(Y) - \phi(Y)|^2 | X=x\right]
		\wedge L')
		\right] \\
		& 
		+ \hat{\mathbb{E}}_S\left[\iint |v(y)-\tilde{v}(x,y)|^2(1 - {\bf 1}_{\mathcal{A}}(y))\nu(dx, dy)\right],
	\end{align*}
	for $L' := 2LQC^2_{p,A}+2C^2_{v, A}$, and $\vartheta(L')$ a function depending on $L'$.
	Note that the set $\mathcal{A}$ is introduced to avoid the restrictive assumption of $v$ being bounded.
	It does not play a role in the actual algorithm.
\end{theorem}
\begin{proof}
	To prepare for our analysis, a more formal construction of the equal partition technique needs to be introduced. 
		
	In practice, for samples $(X^m, Y^m)_{1\leq m \leq M}$ and a measurable sorting function $\mathfrak{S}: \mathbb{R}^q \rightarrow \mathbb{R}$, the $M$ different values can be ordered into, $$\mathfrak{S}(X^{ 1^*}) \leq \mathfrak{S}(X^{2^*}) \leq \cdots \leq \mathfrak{S}(X^{M^*}),$$
	by simply putting $\{X^{1^*}, \cdots, X^{(M/B)^*}\}$ into the first bundle, $\{X^{(M/B+1)^*}, \cdots, X^{(2M/B)^*}\}$ into the second one, etc., assuming $M$ can be divided by $B$ for simplicity.
	
	However, in order to conduct meaningful analysis, a measurable partition of $\mathbb{R}^q$ based on the simulation cloud $(X^m)_{1 \leq m \leq M}$ is required. 
	Thus for any simulation $(X^m)_{1 \leq m \leq M}$ with $\{X^{1^*} = x^*_1, \cdots, X^{M^*} = x^*_M\}$, we define $\mathcal{B}(1) := \mathfrak{S}^{-1}((-\infty, \mathfrak{S}(x^*_{M/B})]), \mathcal{B}(2) := \mathfrak{S}^{-1}((\mathfrak{S}(x^*_{M/B}), \mathfrak{S}(x^*_{2M/B})]), \cdots$, $\mathcal{B}(B) := \mathfrak{S}^{-1}(\mathfrak{S}(x^*_{M- M/B}), \infty)$ and $\cup^B_{b=1}\mathcal{B}(b) = \mathbb{R}^q.$
	Therefore, $$\mathbb{B} = \{\mathfrak{S}^{-1}((-\infty, \mathfrak{S}(x_1)]), \mathfrak{S}^{-1}((\mathfrak{S}(x_1), \mathfrak{S}(x_2)]), \cdots, \mathfrak{S}^{-1}((\mathfrak{S}(x_{B-1}), \infty))\}$$ if and only if 
	\begin{align*}
		&
		(X^{1^{@}}, X^{2^{@}}, \cdots X^{M^{@}}) \\
		\in & 
		(\mathfrak{S}^{-1}((-\infty, \mathfrak{S}(x_1))))^{M/B - 1} \times \{x_1\} \times \cdots \times \{x_{B-1}\} \times (\mathfrak{S}^{-1}((\mathfrak{S}(x_{B-1}), \infty)))^{M/B - 1}.
	\end{align*}
	The notation $@$ denotes any permutation of the set $\{1, 2, \cdots M\}$, noting that each sample is independent of the others and interchangeable. 
	This is measurable with respect to the sigma algebra generated by the simulation cloud $X^m$ as there are finite permutations for fixed $M$ and $\mathfrak{S}$ is measurable. 
		
	Note that this setting is not unique for defining a workable partition and there may be alternative definitions that may improve the analysis result. 
	However, this is an intuitive definition.
		
	Assuming $\sigma(\mathbb{B})$ to be the smallest sigma algebra to determine the partition, we notice that it is smaller than the sigma algebra generated by the random samples $X^m$, $\sigma(\mathbb{B})\subset \sigma(X^m).$
	This is because multiple realizations of the samples can lead to the same partition. 
	A simple thought experiment is to consider a fixed partition, and subsequently move one interior sample within a bundle.
	If we conduct a new bundling with this new set of samples, the partition will remain the same.
	Indeed, the samples within a bundle are independent among each other and have the same distribution.
		
	As for the actual analysis, we start by decomposing the error into different terms for any given partition $\mathbb{B} = \{\mathcal{B}(1), \cdots, \mathcal{B}(B)\} = \{\mathfrak{S}^{-1}((-\infty, \mathfrak{S}(x_1)]), \mathfrak{S}^{-1}((\mathfrak{S}(x_1), \mathfrak{S}(x_2)]), \cdots, \mathfrak{S}^{-1}((\mathfrak{S}(x_{B-1}), \infty))\}$.
	In line with the Monte Carlo literature, we assume $X^i \neq X^j$, if $i \neq j$, and 
	\begin{align}
		&
		\iint |v(y) - \tilde{v}(x,y)|^2 \nu(dx,dy)\nonumber\\
		\leq &
		\sum_{\mathcal{B} \in \mathbb{B}}\int_\mathcal{B}\int|v(y) - \tilde{v}(x, y)|^2 {\bf 1}_{\mathcal{A}}(y)\nu(dx, dy)
		+ \iint |v(y)-\tilde{v}(x,y)|^2(1 - {\bf 1}_{\mathcal{A}}(y))\nu(dx, dy)\nonumber\\
		= &
		\sum^B_{b=1}
		\int_{\mathcal{B}(b)}\int\nu(dx, dy) \bigg(
		||(v - \tilde{v}) {\bf 1}_{\mathcal{A}}||_{\mathcal{B}(b), \infty} 
		-2 ||(v-\tilde{v}){\bf 1}_{\mathcal{A}}||_{\mathcal{B}(b)\backslash \{x_b\}, \#} \nonumber \\
		& \hspace{110pt}
		+2 ||(v-\tilde{v}){\bf 1}_{\mathcal{A}}||_{\mathcal{B}(b)\backslash \{x_b\}, \#}
		-2 ||(v-\tilde{v}){\bf 1}_{\mathcal{A}}||_{\mathcal{B}(b), \#}
		+2 ||(v-\tilde{v}){\bf 1}_{\mathcal{A}}||_{\mathcal{B}(b), \#}
		\bigg)^2\nonumber\\
		&
		+ \iint |v(y)-\tilde{v}(x, y)|^2(1 - {\bf 1}_{\mathcal{A}}(y))\nu(dx, dy)\nonumber\\
		\leq & 
		\sum^B_{b = 1} 
		\int_{\mathcal{B}(b)}\int\nu(dx, dy)  
		\bigg( \max \{||(v-\tilde{v}){\bf 1}_{\mathcal{A}}||_{\mathcal{B}(b), \infty} - 2 ||(v-\tilde{v}){\bf 1}_{\mathcal{A}}||_{\mathcal{B}(b)\backslash\{x_b\}, \#}, 0\} \nonumber\\
		& \hspace{110pt}
		+2 ||(v-\tilde{v}){\bf 1}_{\mathcal{A}}||_{\mathcal{B}(b)\backslash \{x_b\}, \#}
		-2 ||(v-\tilde{v}){\bf 1}_{\mathcal{A}}||_{\mathcal{B}(b), \#}
		+2 ||(v-\tilde{v}){\bf 1}_{\mathcal{A}}||_{\mathcal{B}(b), \#}
		\bigg)^2\nonumber\\
		&
		+ \iint |v(y)-\tilde{v}(x, y)|^2(1 - {\bf 1}_{\mathcal{A}}(y))\nu(dx, dy)\nonumber\\
		\leq & 
		\sum^B_{b = 1}
		\int_{\mathcal{B}(b)}\int\nu(dx, dy)  
		3 \max \{||(v-\tilde{v}){\bf 1}_{\mathcal{A}}||_{\mathcal{B}(b), \infty} - 2 ||(v-\tilde{v}){\bf 1}_{\mathcal{A}}||_{\mathcal{B}(b)\backslash\{x_b\}, \#}, 0\}^2 \nonumber\\
		& +
		\sum^{B-1}_{b = 1}
		\int_{\mathcal{B}(b)}\int\nu(dx, dy)  
		12 (||(v-\tilde{v}){\bf 1}_{\mathcal{A}}||_{\mathcal{B}(b)\backslash \{x_b\}, \#}
		- ||(v-\tilde{v}){\bf 1}_{\mathcal{A}}||_{\mathcal{B}(b), \#})^2 \nonumber\\
		& +
		\sum_{\mathcal{B}\in \mathbb{B}}
		12 \int_\mathcal{B}\int\nu(dx, dy) ||(v-\tilde{v}){\bf 1}_{\mathcal{A}}||^2_{\mathcal{B}, \#}
		+ \iint |v(y)-\tilde{v}(x, y)|^2(1 - {\bf 1}_{\mathcal{A}}(y))\nu(dx, dy)\nonumber\\
		=: & 
		\sum^B_{b = 1} \int_{\mathcal{B}(b)}\int\nu(dx, dy) T_{1, \mathcal{B}(b)}
		+ \sum^{B-1}_{b = 1} \int_{\mathcal{B}(b)}\int\nu(dx, dy) T_{2, \mathcal{B}(b)}\nonumber\\
		&
		+ \sum_{\mathcal{B}\in \mathbb{B}} \int_{\mathcal{B}}\int\nu(dx, dy)  T_{3, \mathcal{B}}
		+ \iint |v(y)-\tilde{v}(x, y)|^2(1 - {\bf 1}_{\mathcal{A}}(y))\nu(dx, dy),
		\label{equation_decomposition}
	\end{align}
	with a slight abuse of notation, $\{x_B\} := \emptyset$ above.
		
	Note that the easiest way to conceptualize the term $||(v-\tilde{v}){\bf 1}_{\mathcal{A}}||_{\mathcal{B}(b)\backslash\{x_b\}, \#}$ is that we simply remove the sample that is used for defining the partition.
	Thus, $$\sum_{m=1}^M{\bf 1}_{\mathcal{B}(b)\backslash\{x_b\}}(X^m) = \sum_{m=1}^M{\bf 1}_{\mathcal{B}(b)}(X^m) - 1.$$
	This is done to ensure that all samples in the empirical norm $||\cdot||_{\mathcal{B}(b)\backslash\{x_b\}, \#}$ are independent of each other.
		
	As the last term cannot be further simplified, we now focus on the first three terms.
	The meaning of all error terms will be discussed in the next subsection.
		
	The first term we study is $T_{3, \mathcal{B}}$, which represents the best possible approximation from the space $H$ to the target function under the empirical norm within the bundle.
	To begin, within any bundle $\mathcal{B}(b)$ under any given partition $\mathbb{B}$, it is obvious that $$||{\bf 1}_{\mathcal{A}}(v - \tilde{v})||_{\mathcal{B}(b), \#} \leq ||v - \tilde{v}_b||_{\mathcal{B}(b), \#} =  \min_{\phi \in H}||v - \phi||_{\mathcal{B}(b), \#},$$ for any $\mathcal{B}(b)$.
	Only the $\tilde{v}_b$ term in the series of $\tilde{v}$ matters and $\tilde{v}_b$ is the function that minimizes the approximation difference under the empirical norm within the given bundle. 
		
	Alternatively, we may consider the following composite norm 
	\begin{equation}
		\sup_{x\in \mathcal{B}(b)}\left(\mathbb{E}\left[ (\mathfrak{f}(Y))^2| X=x\right]\right)^\frac{1}{2}, \label{equation_ad_hoc_norm}
	\end{equation}
	for any given bundle $\mathcal{B}(b)$.
	For the sake of simplicity, we assume there exists an element $\phi_{\mathcal{B}(b)}$ within the space $H$ such that $v(\cdot) - \phi_{\mathcal{B}(b)}(\cdot)$ minimizes the norm, namely, 
	$$\sup_{x\in \mathcal{B}(b)}\mathbb{E}\left[ |v(Y) - \phi_{\mathcal{B}(b)}(Y)|^2 | X=x\right] = \inf_{\phi \in H} \sup_{x\in \mathcal{B}(b)}\mathbb{E}\left[ |v(Y) - \phi(Y)|^2 | X=x\right].$$
	As $\phi_{{\mathcal{B}}(b)} \in H$, it is clear that under the empirical norm, we have $$||{\bf 1}_{\mathcal{A}}(v - \tilde{v})||^2_{\mathcal{B}(b), \#} \leq ||v - \tilde{v}_b||^2_{\mathcal{B}(b), \#}  \leq ||v - \phi_{\mathcal{B}(b)}||^2_{\mathcal{B}(b), \#} = \frac{\sum_{m=1}^{\#\mathcal{B}(b)}|v(Y^{b, m}) - \phi_{\mathcal{B}}(Y^{b, m})|^2}{\# \mathcal{B}(b)}.$$
		
	Without loss of generality, assume $(X^{b, \#\mathcal{B}(b)}, Y^{b, \#\mathcal{B}(b)})$ is the bundle defining sample as stated in the construction of equal partition if $b \neq B$. 
	Recalling that samples within a bundle are i.i.d. given the partition, we can take the conditional expectation of the empirical norm with respect to the position of $(X^{b, \#\mathcal{B}(b)}, Y^{b, \#\mathcal{B}(b)})_{1\leq m \leq \#\mathcal{B}(b) - 1}$.
	\begin{align*}
		&
		\hat{\mathbb{E}}\left[\left.||{\bf 1}_{\mathcal{A}}(v - \tilde{v})||^2_{\mathcal{B}, \#}\right|\sigma(\mathbb{B}), X^{b, 1}, X^{b, 2}, \cdots, X^{b, \#\mathcal{B}(b) -1}\right] \\ 
		\leq &
		\frac{\sum_{m=1}^{\#\mathcal{B}(b)}\hat{\mathbb{E}}\left[\left.|v(Y^{b, m}) - \phi_{\mathcal{B}}(Y^{b, m})|^2\right|\sigma(\mathbb{B}),  X^{b, 1}, \cdots, X^{b, \#\mathcal{B}(b)-1}\right]}{\# \mathcal{B}(b)}\\
		= &
		\frac{\sum_{m=1}^{\#\mathcal{B}(b) - 1}\mathbb{E}\left[\left.|v(Y) - \phi_{\mathcal{B}}(Y)|^2\right| X = X^{b, m} \in \mathcal{B}(b)\backslash\{x_b\}\right] + \mathbb{E}\left[\left.|v(Y) - \phi_{\mathcal{B}}(Y)|^2\right| X = x_b\right]}{\# \mathcal{B}(b)}\\
		\leq &
		\frac{\sum_{m=1}^{\#\mathcal{B}(b)}\sup_{x \in \mathcal{B}(b)}\mathbb{E}\left[\left.|v(Y) - \phi_{\mathcal{B}}(Y)|^2\right| X =x \right]}{\# \mathcal{B}(b)}\\
		= &
		\inf_{\phi \in H}\sup_{x \in \mathcal{B}(b)}\mathbb{E}\left[\left.|v(Y) - \phi(Y)|^2\right| X = x \right].
	\end{align*}
	There are some details in the above calculation that require explanation.
	Note that the boundary point information is included in $\sigma(\mathbb{B})$, therefore we only calculate conditional expectations with the remaining samples. 
	For the last bundle $\mathcal{B}(B)$, the separation is not necessary as no sample is used to define the bundle, but this does not alter the result. 
	We use the fact that each sample is independent within the bundle in the first equality. 
		
	Next, if the minimal element $\phi_{\mathcal{B}(b)}$ does not exist, one has to adjust the derivation above with a limiting argument. 
	This means that by the definition of infimum, one can find a sequence of functions $(v - \phi_{\mathcal{B}(b) , n})_{n \in \mathbb{Z}^+}$, such that 
	$$\sup_{x\in \mathcal{B}(b)}\mathbb{E}\left[ |v(Y) - \phi_{\mathcal{B}(b),n}(Y)|^2 | X=x\right] \leq \inf_{\phi \in H} \sup_{x\in \mathcal{B}(b)}\mathbb{E}\left[ |v(Y) - \phi(Y)|^2 | X=x\right] + \frac{1}{n}.$$
	By repeating the above argument for each function in the sequence, replacing the infimum in the proof by the corresponding upper bound and taking $n$ to infinity (with the eventual inequality), we arrive at the same conclusion.
		
	Thereafter, if we consider that expectation of the empirical norm conditioning on $\sigma(\mathbb{B})$, we have 
	\begin{align*}
		\hat{\mathbb{E}}\left[\left.||{\bf 1}_{\mathcal{A}}(v - \tilde{v})||^2_{\mathcal{B}, \#}\right|\sigma(\mathbb{B})\right]
		= & 
		\hat{\mathbb{E}}\left[\left.|\hat{\mathbb{E}}\left[\left.||{\bf 1}_{\mathcal{A}}(v - \tilde{v})||^2_{\mathcal{B}, \#}\right|\sigma(\mathbb{B}), X^{b, 1}, X^{b, 2}, \cdots, X^{b, \#\mathcal{B}(b) -1}\right] \right|\sigma(\mathbb{B})\right]\\ 
		\leq &
		\inf_{\phi \in H}\sup_{x \in \mathcal{B}(b)}\mathbb{E}\left[\left.|v(Y) - \phi(Y)|^2\right| X = x \right].
	\end{align*}
		
	Note that this bound is defined on all given partitions and solely depends on the partition but not the choice of $\phi_{\mathcal{B}}$ for any bundle $\mathcal{B}$.
	Our calculation here is purely within a bundle given the partition is known. 
	Therefore, even if the minimum function $\phi_{\mathcal{B}(b)}$ or the sequence $(\phi_{\mathcal{B}(b) , n})_{n \in \mathbb{Z}^+}$ is not unique, the actual choice of these functions does not matter as long as they are picked in a consistent and measurable way. 
		
	Here, we derive a bound for the expectation of the weight summation $T_{3,\mathcal{B}}$ in \eqref{equation_decomposition} with respect to the simulation cloud, i.e.,
	\begin{align}
		\hat{\mathbb{E}}_S\left[\sum_{\mathcal{B}\in \mathbb{B}} \int_{\mathcal{B}}\int\nu(dx, dy) T_{3,\mathcal{B}}\right] 
		& \leq 
		12 \hat{\mathbb{E}}_S\left[
		\sum_{\mathcal{B}\in \mathbb{B}}\int_\mathcal{B}\int\nu(dx,dy)||v - \phi_\mathcal{B}||^2_{\mathcal{B}, \#}
		\right]\nonumber\\
		& \leq 
		\frac{12}{\hat{\mathbb{E}}[{\bf 1}_S]}\hat{\mathbb{E}}
		\left[
		\sum_{\mathcal{B}\in \mathbb{B}}\int_\mathcal{B}\int\nu(dx,dy)					\hat{\mathbb{E}}\left[\left.||v - \phi_\mathcal{B}||^2_{\mathcal{B}, \#}\right|\sigma(\mathbb{B})\right]
		\right]\label{equation_T2}\\
		& \leq  
		\frac{12}{\hat{\mathbb{E}}[{\bf 1}_S]}
		\hat{\mathbb{E}}\left[
		\sum_{\mathcal{B}\in \mathbb{B}}\int_\mathcal{B}\int\nu(dx,dy)\inf_{\phi \in H} \sup_{x\in \mathcal{B}}\mathbb{E}\left[ |v(Y) - \phi(Y)|^2 | X=x\right]
		\right].\nonumber
	\end{align}
	In this inequality, we expand the denominator of our adjusted probability by also including the rejected cases and applying the results above for each partition.
	
	However, for an unbounded bundle $\mathcal{B}$, it is possible to find an example such that \\ $\sup_{x\in \mathcal{B}}\mathbb{E}\left[ |v(Y) - \phi(Y)|^2 | X=x\right] = \infty $.
	An alternative bound is required to ensure that our error bound is not trivial.
	Note that given the square norm $|\alpha(b)|^2 \leq L$, we have $$\forall y, b, \,|\tilde{v}_{\mathcal{B}(b)}(y){\bf 1}_{\mathcal{A}}(y)|^2 \leq \left(\sum^Q_{l=1}|\alpha_l(b)|^2\right)\left(\sum^Q_{l=1}|p_l(y){\bf 1}_{\mathcal{A}}(y)|^2\right)\leq L Q\max_{l =1 \ldots, Q}\max_{y\in \mathcal{A}}|p_l(y)|^2,$$
	and
	\begin{align*}
		||(v-\tilde{v}){\bf 1}_{\mathcal{A}}||^2_{\mathcal{B}, \#}
		&
		=\frac{\sum_{m=1}^M{\bf 1}_{\mathcal{B}}(X^m)|(v(Y^m)-\phi_\mathcal{B}(Y^m)){\bf 1}_{\mathcal{A}}|^2}{\sum_{m=1}^M{\bf 1}_{\mathcal{B}}(X^m)}
		\nonumber\\
		&
		\leq \frac{\sum_{m=1}^M{\bf 1}_{\mathcal{B}}(X^m)(2|(v(Y^m){\bf 1}_{\mathcal{A}}|^2 + 2|\phi_\mathcal{B}(Y^m)){\bf 1}_{\mathcal{A}}|^2}{\sum_{m=1}^M{\bf 1}_{\mathcal{B}}(X^m)}\nonumber\\
		&
		\leq \frac{\sum_{m=1}^M{\bf 1}_{\mathcal{B}}(X^m)(2C_{v,\mathcal{A}}^2 + 2LQC_{p, \mathcal{A}}^2)}{\sum_{m=1}^M{\bf 1}_{\mathcal{B}}(X^m)} = L'.
	\end{align*}
	We shall use this alternative bound in Equation \eqref{equation_T2} for any bundle $\mathcal{B}$ such that \\ $\sup_{x\in \mathcal{B}}\mathbb{E}\left[ |v(Y) - \phi(Y)|^2 | X=x\right] > L'$.
	Combining the two error bounds, we have 
	\begin{align}
		\hat{\mathbb{E}}_S\left[\sum_{\mathcal{B} \in \mathbb{B}}T_{3,\mathcal{B}}\right] 
		& \leq  
		\frac{12}{\hat{\mathbb{E}}[{\bf 1}_S]}
		\hat{\mathbb{E}}\left[
		\sum_{\mathcal{B}\in \mathbb{B}}\int_\mathcal{B}\int\nu(dx,dy)
		(\inf_{\phi \in H} \sup_{x\in \mathcal{B}}\mathbb{E}\left[ |v(Y) - \phi(Y)|^2 | X=x\right]
		\wedge L')
		\right].\nonumber
	\end{align}
		
	Next, we consider the term $T_{1,\mathcal{B}}$ in \eqref{equation_decomposition}.
	This term concerns the difference between the theoretical projection and the empirical regression function within each bundle. 
	Here we restate that $S$ denotes the modified probability, based on the regression coefficients, where $\mathcal{A}$ is a compact set defined with respect to $Y$ only. 
		
	By taking conditional expectations with respect to $\sigma(\mathbb{B})$, we have,
	\begin{align*}
		&
		\hat{\mathbb{E}}_S\left[\sum^B_{b = 1}\int_{\mathcal{B}(b)}\int\nu(dx, dy)T_{1,\mathcal{B}(b)}\right] \\
		= &
		\frac{1}{\hat{\mathbb{E}}[{\bf 1}_S]}\hat{\mathbb{E}}\left[
		\sum^B_{b = 1}\int_{\mathcal{B}(b)}\int\nu(dx, dy)
		\hat{\mathbb{E}}[
		T_{1,\mathcal{B}(b)}{\bf 1}_S
		|\sigma(\mathbb{B})]
		\right]\\
		= & 
		\frac{1}{\hat{\mathbb{E}}[{\bf 1}_S]}\hat{\mathbb{E}}\left[
		\sum^B_{b = 1}\int_{\mathcal{B}(b)}\int\nu(dx, dy)
		\hat{\mathbb{E}}_{\mathcal{B}(b)}[
		{\bf 1}_S 3 \max \{||{\bf 1}_{\mathcal{A}}(v-\tilde{v})||_{\mathcal{B}(b), \infty} 
		- 2 ||(v-\tilde{v}){\bf 1}_{\mathcal{A}}||_{\mathcal{B}(b)\backslash\{x_b\}, \#}, 0\}^2
		]
		\right].
	\end{align*}
	It is important that we condition on the smaller sigma algebra such that all samples have an identical conditional distribution.
	If we can condition on the whole sigma algebra generated by $(X^m)_{1\leq m \leq M}$, each $Y^m$ will have a different distribution depending on the position of $X^m$.
	Within each (given) bundle $\mathcal{B}(b)$, we may consider the two norms $||\cdot||_{\mathcal{B}(b),\infty}$ and $||\cdot||_{\mathcal{B}(b)\backslash\{x_b\}, \#}$ as the theoretical and empirical $L^2$ norms of a random process satisfying the probability distribution $\mathbb{P}_{\mathcal{B}(b)}:= \frac{\int_{\mathcal{B}(b)}\nu(dx, \cdot)}{\int_{\mathcal{B}(b)}\int\nu(dx, dy)}$ and extend our notation for expectations to this measure. 
	In other words, $\mathbb{P}_{\mathcal{B}(b)}$ is the conditional probability of $Y^m$ given that $X^m$ is within the bundle.
	As only the samples within bundle $\mathcal{B}(b)$ are considered in $T_{1, \mathcal{B}(b)}$, we only have to consider the identically distributed samples following $\mathbb{P}_{\mathcal{B}(b)}$.
	Thus, we simplify the notation with this measure.
		
	Assume that $\sum^M_{m=1}{\bf 1}_{\mathcal{B}(b)\backslash\{x_b\}}(X^m) = N - 1$ and let $u > 864 L'/(N - 1)$ be arbitrary, by Theorem 11.2 in \cite{gyorfi_kohler_krzyzak_walk_2002}, we find
	\begin{align}
		& 
		\hat{\mathbb{P}}_{\mathcal{B}(b)}\{3 \max 
		\{||{\bf 1}_{\mathcal{A}}(v-\tilde{v})||_{\mathcal{B}(b), \infty} 
		- 2 ||(v-\tilde{v}){\bf 1}_{\mathcal{A}}||_{\mathcal{B}(b)\backslash\{x_b\}, \#}, 0\}^2 > u 
		\mbox{ and the event } S \mbox{ is true}\}\nonumber\\
		\leq & 
		\hat{\mathbb{P}}_{\mathcal{B}(b)}\{ \exists \phi \in H_L : 
		||{\bf 1}_{\mathcal{A}}(v-\phi)||_{\mathcal{B}(b), \infty} 
		- 2 ||(v-\phi){\bf 1}_{\mathcal{A}}||_{\mathcal{B}(b)\backslash\{x_b\}, \#} > \sqrt{u/3} \mbox{ and the event } S \mbox{ is true}\}\nonumber\\
		\leq & 
		\hat{\mathbb{P}}_{\mathcal{B}(b)}\{ \exists \phi \in H_L : 
		||{\bf 1}_{\mathcal{A}}(v-\phi)||_{\mathcal{B}(b), \infty} 
		- 2 ||(v-\phi){\bf 1}_{\mathcal{A}}||_{\mathcal{B}(b)\backslash\{x_b\}, \#} > \sqrt{u/3} \}\nonumber\\
		\leq &
		3\hat{\mathbb{E}}_{\mathcal{B}(b)}[\mathcal{N}_2(\sqrt{2/3}\sqrt{u}/24, H_{L,\mathcal{A}}, Y_{\mathcal{B}(b)}^{2(N-1)})]\exp\left(-\frac{(N-1)u}{864L'}\right)\nonumber\\
		\leq &
		3\hat{\mathbb{E}}_{\mathcal{B}(b)}[\mathcal{N}_2(\sqrt{L'}/\sqrt{N-1}, H_{L,\mathcal{A}}, Y_{\mathcal{B}(b)}^{2(N-1)})]\exp\left(-\frac{(N-1)u}{864L'}\right), \label{equation_covering_number}
	\end{align}
	where $H_L$ is the set of all functions in $H$ whose coordinates with respect to the basis $(p_l)_{1\leq l \leq Q}$, have a Euclidean norm no greater than $L$ and $H_{L, \mathcal{A}}$ the set containing all functions of the form ${\bf 1}_{\mathcal{A}}(\phi - v)$, where $\phi$ belongs to $H_L$.
	Again, since we condition only on the partition, all samples within a bundle are i.i.d. and the condition for Theorem 11.2 in \cite{gyorfi_kohler_krzyzak_walk_2002} is satisfied. 
	In fact, this proof should work for all partitions for which the samples remain i.i.d. within a bundle.
	Note that the indicator for the event $S$ is kept in the first line to keep our regression function bounded, then we drop the indicator in the second inequality to take advantage of the independent samples.
	Finally, $Y_{\mathcal{B}(b)}$ is a sample following the conditional probability $\mathbb{P}_\mathcal{B}$ here.
		
	Constant $\mathcal{N}_2$ in \eqref{equation_covering_number} is called the covering number and it is bounded by Lemma 9.2 and Theorem 9.4 of \cite{gyorfi_kohler_krzyzak_walk_2002}: $$\mathcal{N}_2(\sqrt{L'}/\sqrt{N-1}, H_{L,\mathcal{A}}, Y^{2(N-1)}) \leq 3\left(\frac{2eL'}{L'/(N-1)}\log(\frac{3eL'}{L'/(N-1)})\right)^{V_{H^+_{L, \mathcal{A}}}}\leq 3\left[3e(N-1)\right]^{2V_{H^+_{L, \mathcal{A}}}},$$ where $V$ denotes the Vapnik-Chervonenkis dimension, which represents the number of elements in the largest set that can be shattered by a class of subsets in $\mathbb{R}^q$.
	The reader is referred to section 9.4 of \cite{gyorfi_kohler_krzyzak_walk_2002} for further information on $\mathcal{N}_2$ and $V$.
		
	Next, recalling the definition $H^+$ from Section \ref{section_introduction}, we notice that $V_{H^+_{L, \mathcal{A}}} \leq V_{H^+_L}$, which can be shown by the following argument.
	Let $(y, z) \in \mathbb{R}^{q} \times \mathbb{R}$, if $y \not\in \mathcal{A}$ and $z\geq 0$, then $(y,z)$ is contained in none of the sets in $H_{L,\mathcal{A}}^+$ and if $y \not\in \mathcal{A}$ and $z \leq 0$, then $(y,z)$ is contained in each set of $H^+_{L,\mathcal{A}}$.
	Hence, if $H^+_{L,\mathcal{A}}$ shatters a set of points, then the x-coordinates of these points must lie in $\mathcal{A}$ and $H^+_L$ also shatters this set of points. 
	
	In addition, we have the fact that $H_L \subset H$ and observe that 
	\begin{align*}
	H^+ 
	& \subseteq \{\{(x,t) :  \phi(x) + a_0t \geq 0\} : \phi\in H, a_0 \in \mathbb{R}\},
	\end{align*}
	which is a linear vector space of dimension less than or equal to $Q+1$, thus Theorem 9.5 of \cite{gyorfi_kohler_krzyzak_walk_2002} implies $$V_{H^+_L}\leq Q+1.$$
	It follows that, for any $u > 864L'/(N-1)$, the probability under consideration is bounded by \\ $9 [3e(N-1)]^{2(Q+1)}\exp\left(-\frac{(N-1) u}{864L'}\right),$
	and for any $w>864L'/(N-1)$,
	\begin{align*}
		&
		\hat{\mathbb{E}}[T_{1,\mathcal{B}(b)}{\bf 1}_S|\sigma(\mathbb{B}), \sum^M_{m=1}{\bf 1}_{\mathcal{B}(b)\backslash\{x_b\}}(X^m) = N-1] \\
		\leq & 
		w + 9[3e(N-1)]^{2(Q+1)} \int^\infty_w\exp\left(-\frac{(N-1)t}{864L'}\right)dt \\
		& =
		w + 9[3e(N-1)]^{2(Q+1)}\frac{864 L'}{N-1}\exp\left(-\frac{(N-1)w}{864L'}\right).
	\end{align*}
	By setting, $$w = \frac{864L'}{N-1} \log\left(9[3e(N-1)]^{2(Q+1)}\right),$$ and taking expectations with respect to $\sigma(\mathbb{B})$, we find
	\begin{align*}
		&
		\hat{\mathbb{E}}_S\left[\sum^B_{b =1}\int_{\mathcal{B}(b)}\int\nu(dx, dy)T_{1,\mathcal{B}(b)}\right]\\
		\leq & 
		\frac{\vartheta(L')}{\hat{\mathbb{E}}[{\bf 1}_S]}\hat{\mathbb{E}}
		\left[\sum^B_{b = 1}\int_{\mathcal{B}(b)}\int\nu(dx, dy)\frac{(\log(\sum^M_{m=1}{\bf 1}_{\mathcal{B}(b)}(X^m) - 1)+1)(Q+1)}{\sum^M_{m=1}{\bf 1}_{\mathcal{B}(b)}(X^m) - 1}\right],
		\end{align*}
	where one possible choice of $\vartheta$ is $\vartheta(L'):=1728(\log(27e) + 1) L'$.
	This can be checked by simple algebra.
	Note that $\vartheta$ is independent of the number of samples in a bundle and only depends on $L'$.
	
	Finally, $T_{2, \mathcal{B}(b)}$ is the technical term introduced by the definition of the partition $\mathbb{B}$.
	Consider  any realization of the simulation cloud $(X^m, Y^m)$ and partition $\mathbb{B}$, in particular, there exist boundary defining samples $(x_b, y_b) \in (X^m, Y^m)$ for $b = 1, 2, \cdots, B-1$.
	Using the inequality $(\sqrt{a} - \sqrt{b})^2 \leq |a-b|$ and the definition of $L'$, we have
	\begin{align*}
		&
		12 (||(v-\tilde{v}){\bf 1}_{\mathcal{A}}||_{\mathcal{B}(b)\backslash \{x_b\}, \#}
		- ||(v-\tilde{v}){\bf 1}_{\mathcal{A}}||_{\mathcal{B}(b), \#})^2 \\
		\leq & 
		12 \left|\frac{\sum_{m=1}^M{\bf 1}_{\mathcal{B}(b)\backslash\{x_b\}}(X^m)|(v(Y^m)-\tilde{v}(X^m, Y^m){\bf 1}_\mathcal{A}(Y^m)|^2}{\sum_{m=1}^M{\bf 1}_{\mathcal{B}(b)}(X^m) - 1}\right.\\
		& \hspace{20pt} \left.
		- \frac{\sum_{m=1}^M{\bf 1}_{\mathcal{B}(b)}(X^m)|(v(Y^m)-\tilde{v}(X^m, Y^m){\bf 1}_\mathcal{A}(Y^m)|^2}{\sum_{m=1}^M{\bf 1}_{\mathcal{B}(b)}(X^m)}\right|\\
		\leq &
		12 \left|\frac{\sum_{m=1}^M{\bf 1}_{\mathcal{B}(b)\backslash\{x_b\}}(X^m)(|(v(Y^m)-\tilde{v}(X^m, Y^m){\bf 1}_\mathcal{A}(Y^m)|^2 - |(v(y_m)-\tilde{v}(x_m, y_m){\bf 1}_\mathcal{A}(y_m)|^2)}
		{(\sum_{m=1}^M{\bf 1}_{\mathcal{B}(b)}(X^m) - 1)(\sum_{m=1}^M{\bf 1}_{\mathcal{B}(b)}(X^m))}\right|\\
		\leq &
		\frac{24 L'}{(\sum_{m=1}^M{\bf 1}_{\mathcal{B}(b)}(X^m))},
	\end{align*}
	for $b = 1, 2, \cdots, B-1.$
		
	Substituting all the results above into Equation \eqref{equation_decomposition} we conclude the proof.
\end{proof}

\begin{remark} \label{remark_empty_bundle_equal}
	Implicitly, it is assumed in our proof that the sorting function used behaves nicely such that there is no empty bundle or bundle with a few samples in our partition. 
	If this is not the case, these bundles need to be merged with other bundles in a consistent way and the number of bundles needs to be adjusted accordingly. 
	We omit this extra complexity in favor of the presentation.
\end{remark}

\subsection{Discussion on the Error Bound}
We shall discuss the meaning of all error terms in Theorem \ref{theorem_regression}.
Note that most discussions here are in heuristic sense instead of rigorous  analysis.

The last term in the sum $\hat{\mathbb{E}}_S\left[\iint |v(y)-\tilde{v}(x,y)|^2(1 - {\bf 1}_{\mathcal{A}}(y))\nu(dx, dy)\right]$ represents an expectation with respect to an integration of the approximation error using the probability measure of $(X, Y)$ outside a given compact set $\mathcal{A}$ for $Y$.
In theory, for an increasing series of compact sets $\mathcal{A}_1 \subset \mathcal{A}_2 \subset \ldots \subset \mathcal{A}_A \subset \mathbb{R}^q$, $$\lim_{A\rightarrow \infty}\iint |v(y)-\tilde{v}(x, y)|^2(1 - {\bf 1}_{\mathcal{A}_A}(y))\nu(dx, dy)\rightarrow 0,$$ since the original function and the approximant are both in $L^2(\nu)$ and the dominating convergent theorem. 
Again, the set $\mathcal{A}$ plays no role in the algorithm. 
This term is introduced to reflect that only the area with high probability measure has strong impact on a Monte Carlo approximation, thus, we can only consider the behavior of function $v$ in this area and do not impose strong conditions on $v$ over the whole domain.   
In practice, we can find a big enough set $\mathcal{A}$ such that the last term is smaller than any preset tolerated level.
Otherwise, the probability distribution is too much spread out so that Monte-Carlo may not be a suitable approximation method.

The term $\frac{12}{\hat{\mathbb{E}}[{\bf 1}_S]}
	\hat{\mathbb{E}}\left[
	\sum_{\mathcal{B}\in \mathbb{B}}\int_\mathcal{B}\int\nu(dx,dy)
	(\inf_{\phi \in H} \sup_{x\in \mathcal{B}}\mathbb{E}\left[ |v(Y) - \phi(Y)|^2 | X=x\right]
	\wedge L') \right]$ can be seen as the average of the best projection error among bundles and upper bounded by $L'$. 
It concerns how close the original function and its projection onto the space spanned by our basis are.
This term should be controlled by increasing the number of bundles. 
As the number of bundles increases, all the bundles converge to a point and the error term becomes $\frac{12}{\hat{\mathbb{E}}[{\bf 1}_S]}
	\hat{\mathbb{E}}\left[
	(\inf_{\phi \in H} \mathbb{E}\left[ |v(Y) - \phi(Y)|^2 | X\right]
	\wedge L')
	\right]$, which is the average projection error over the whole range of $X$. 
It is clear that $\hat{\mathbb{E}}[{\bf 1}_S]$ might change when we increase the number of bundles, thus the above analysis is by no mean rigorous.
We will provide comments on $S$ at the end of this section.

The first error term is the bound for the estimation error based on the empirical norm instead of the theoretical norm.
This term can be controlled by simply increasing the number of samples within each bundle such that this term is below a certain threshold.
Because low sample numbers imply high error bounds, if there are bundles in any partition which contain very few samples, they should be merged with other bundles or the sorting function must be adapted. 
	
The second error term is just a technical term for constructing a measurable partition and can be controlled by the number of samples in each bundle.
	
Therefore, the best way to set the parameters for the SGBM algorithm is to first fix the number of samples in each bundle such that the first two terms in the error bound are under a given threshold, then increase the number of bundles to control the projection error.
However, we cannot write a simple convergence rate for the combined error under the current conditions.

\subsection{Discussion on Event $S$}
There is one final question remaining. 
The above argument heavily depends on $L$, which is a user defined quantity, and the probability of $S$, the event that $L^2$ norms of the regression coefficients are below threshold $L$.
It is natural to ask if we can actually find a number $L$ such that $\mathbb{P}_S$ is bounded below, as our bound becomes trivial when $\frac{1}{\hat{\mathbb{E}}[{\bf 1}_S]}$ tends to infinity and it would be incredibly expensive to apply the algorithm if we reject most of the simulations.
In the following, we shall provide a heuristic argument for the convergence of the regression coefficients within the bundle, therefore, there is a natural choice of $L$ depending on the target function itself and a hard cut off may be unnecessary.
The numerical experiments in the next section also back up this argument. 

Once again, we consider Equation \eqref{equation_coefficients}, where the regression coefficients within any bundle $\mathcal{B}(b)$ satisfy
\begin{equation*}
\mathcal{I}^\top\mathcal{I}\alpha(b) = \mathcal{I}^\top v(Y^b),
\end{equation*}
with 
\begin{equation*}
\mathcal{I}^\top\mathcal{I} = \left(
\begin{array}{ccc}
\sum^{\# \mathcal{B}(b)}_{i = 1}(p_1(Y^{b, i}))^2 & & \sum^{\# \mathcal{B}(b)}_{i = 1}p_1(Y^{b, i})p_Q(Y^{b, i}) \\
& \ddots &\\
\sum^{\# \mathcal{B}(b)}_{i = 1}p_Q(Y^{b, i})p_1(Y^{b, i}) & & \sum^{\# \mathcal{B}(b)}_{i = 1}(p_Q(Y^{b, i}))^2 \\
\end{array} \right)
\end{equation*}
and 
\begin{equation*}
\mathcal{I}^\top v(Y^b) = \left(\begin{array}{c}
\sum^{\# \mathcal{B}(b)}_{i = 1}p_1(Y^{b, i})v(Y^{b, i})\\
\vdots\\
\sum^{\# \mathcal{B}(b)}_{i = 1}p_Q(Y^{b, i})v(Y^{b, i})
\end{array}\right).
\end{equation*}
When the number of samples within a bundle tends to infinity, it is easy to see that 
\begin{equation*}
\frac{1}{\#\mathcal{B}(b)}\mathcal{I}^\top\mathcal{I} \rightarrow \left(
\begin{array}{ccc}
\mathbb{E}[p_1(Y))^2|X\in \mathcal{B}(b)] & & \mathbb{E}[p_1(Y)p_Q(Y)|X\in\mathcal{B}(b)] \\
& \ddots &\\
\mathbb{E}[p_Q(Y)p_1(Y)|X\in\mathcal{B}(b)] & & \mathbb{E}[p_Q(Y))^2|X\in \mathcal{B}(b)]\\
\end{array} \right)
\end{equation*}
and 
\begin{equation*}
\frac{1}{\#\mathcal{B}(b)}\mathcal{I}^\top v(Y^b) \rightarrow \left(\begin{array}{c}
\mathbb{E}[p_1(Y)v(Y)|X\in\mathcal{B}(b)]\\
\vdots\\
\mathbb{E}[p_Q(Y)v(Y)|X\in\mathcal{B}(b)]
\end{array}\right).
\end{equation*}
So, the empirical system of equations "converges" to the system of equations of a projection.
Therefore, as long as our basis is properly defined such that they remain linearly independent for all bundles, this system of equations should be solvable with enough simulation paths.
Moreover, since the regression coefficients should "converge" to the theoretical projection coefficients, we could pick $L$ depending on the $L^2$ norm of $v$ itself, like, for example, two times its theoretical norm.
Alternatively, when there are enough samples within each bundle, we suspect that the regression coefficients simply "converges" to the theoretical value and satisfy the bounded condition of regression in a natural way.
Therefore, no actual rejection step in the algorithm is needed when there are enough samples within the bundles.
This proposition appears to be supported by our numerical experiments.

However, there are multiple difficulties to incorporate the above argument into Theorem \ref{theorem_regression}. 
First, since equal partitioning is not a recursive partitioning scheme as defined in \cite{gordon1984almost}, we cannot use a martingale argument on equal partitioning, limiting the available tools.
Secondly, as the partition changes when increasing the overall number of samples, the above "convergence" does not seem to be properly defined.
Finally, we would have to introduce a measure of convergence with respect to a matrix inverse, which is beyond the scope of this work. 

On the other hand, there is a possibility that when the number of samples within a bundle is too small, the algorithm as a whole will fail to converge.
Thus, it is beneficial to remind a user of such possibility and put a safety check in place.
Therefore, we keep the derivation of Theorem \ref{theorem_regression} as a complete justification for SGBM. In practice, one can either make sure that there are enough samples within each bundle and let $L$ be infinity.
In this case, Theorem \ref{theorem_regression} no longer applies but we believe that the overall error will satisfy a bound of similar form. 
Alternatively, one starts from a small $L$ when running the algorithm and increases $L$'s value until most tests are accepted.
By these techniques, the error bound from Theorem \ref{theorem_regression} remains valid.

\section{Error Analysis} \label{section_explicit}

A complete error description of the algorithm with respect to the application of SGBM towards BSDEs will be derived in this section.

We  wish to apply the theorem from the last section to establish an error bound for the expectation of our approximation with respect to the selected simulation cloud.
We need to check that after rejecting the simulations that generate regression coefficients that are "too large", our approximation functions are bounded in the recursion.
We notice that for any $k\leq N$,
$$|y^{R}_k(x)| \leq \max \{C_{M, A} L \sqrt{2(1 + C_\pi^2)},  C_{\Phi, A} \}=: C_{Y,A}$$ and  $$|z^{R}_k(x)| \leq C_{M, A} L \sqrt{2(1 + C_\pi^2)} =: C_{Z, A}$$ for all $x$ in a compact set $A$.
The constant $C_\pi$ is defined as $\max_{k = 0, \ldots, N-1}\Delta_k$. 
These bounds can be proven by Assumption \ref{assumption:globally_lipschitz} and some simple inequalities.
Furthermore, we have $\forall x \in A$, $$f^{R}_k(x) := f_{k}(y^{R}_{k}(x), z^{R}_{k}(x)) \leq C_f + L_f (C_{Y, A} + C_{Z, A}) =: C_{f,A},$$ which follows from the Lipschitz assumptions of $f$.
Therefore, Theorem \ref{theorem_regression} applies.

We denote by $S$ the set of all simulation cloud values $(X^{\pi,m}_{t_k})_{\substack{1\leq m \leq M\\ 0 \leq k \leq N}}$ such that the Euclidean norm of the regression coefficients at each time step in each bundle is bounded by $L$, and the expectation is adjusted accordingly. 
With the application of Theorem \ref{theorem_regression}, we know that for any given compact set $\mathcal{A}$,
\begin{align*}
	&
	\hat{\mathbb{E}}^x_{t_{k}, S}\left[\expectation{|y^{R}_{i+1}(X^\pi_{t_{i+1}}) - \tilde{y}^{R}_{i+1}(X^\pi_{t_i}, X^\pi_{t_{i+1}})|^2}{k}{x}\right]\\
	\leq &
	\frac{\vartheta(L'_y)}{\hat{\mathbb{E}}^x_{t_k}[{\bf 1}_S]}
	\hat{\mathbb{E}}^x_{t_k}\left[\sum^B_{b = 1}\int_{\mathcal{B}(b)}\int\nu(dx, dy)\frac{(\log(\sum^M_{m=1}{\bf 1}_{\mathcal{B}(b)}(X^m) - 1)+1)(Q+1)}{\sum^M_{m=1}{\bf 1}_{\mathcal{B}(b)}(X^m) - 1}\right]\\
	& +
	\hat{\mathbb{E}}^x_{t_k, S}\left[\sum^{B-1}_{b =1}\int_{\mathcal{B}(b)}\int\nu(dx, dy)\frac{24 L'_y}{(\sum_{m=1}^M{\bf 1}_{\mathcal{B}(b)}(X^m))}\right]\\
	& + 
	\frac{12}{\hat{\mathbb{E}}^x_{t_k}[{\bf 1}_S]}
	\hat{\mathbb{E}}^x_{t_k}\left[
	\sum_{\mathcal{B}\in \mathbb{B}}\int_\mathcal{B}\int\nu(dx,dy)
	(\inf_{\phi \in H} \sup_{\theta\in \mathcal{B}}\mathbb{E}^\theta_{t_i}\left[|y^{R}_{t_{i+1}}(X^\pi_{t_{i+1}}) - \phi(X^\pi_{t_{i+1}})|^2\right]
	\wedge L'_y)
	\right] \\
	& + \hat{\mathbb{E}}^x_{t_k, S}\left[\expectation{|y^{R}_{i+1}(X^\pi_{t_{i+1}}) - \tilde{y}^{R}_{i+1}(X^\pi_{t_i}, X^\pi_{t_{i+1}})|^2(1 - {\bf 1}_{\mathcal{A}}(X^\pi_{t_{i+1}}))}{k}{x}\right] =: \Xi^x_{t_k}(i,y).
\end{align*}
and 
\begin{align*}
	& 
	\hat{\mathbb{E}}^x_{t_k, S}\left[\expectation{
		|f^{R}_{i+1}(X^\pi_{t_{i+1}})
		- \tilde{f}^{R}_{i+1}(X^\pi_{t_i}, X^\pi_{t_{i+1}})|^2
	}{k}{x}\right]\\
	\leq &
	\frac{\vartheta({L'_f})}{\hat{\mathbb{E}}^x_{t_k}[{\bf 1}_S]}\hat{\mathbb{E}}^x_{t_k}
	\left[\sum^B_{b = 1}\int_{\mathcal{B}(b)}\int\nu(dx, dy)\frac{(\log(\sum^M_{m=1}{\bf 1}_{\mathcal{B}(b)}(X^m) - 1)+1)(Q+1)}{\sum^M_{m=1}{\bf 1}_{\mathcal{B}(b)}(X^m) - 1}\right] \\
	& +
	\hat{\mathbb{E}}^x_{t_k, S}\left[\sum^{B-1}_{b =1}\int_{\mathcal{B}(b)}\int\nu(dx, dy)\frac{24 L'_f}{(\sum_{m=1}^M{\bf 1}_{\mathcal{B}(b)}(X^m))}\right]\\
	& +
	\frac{12}{\hat{\mathbb{E}}^x_{t_k}[{\bf 1}_S]}
	\hat{\mathbb{E}}^x_{t_k}\left[
	\sum_{\mathcal{B}\in \mathbb{B}}\int_\mathcal{B}\int\nu(dx,dy)
	(\inf_{\phi \in H} \sup_{\theta\in \mathcal{B}}\mathbb{E}^\theta_{t_i}\left[|f^{R}_{t_{i+1}}(X^\pi_{t_{i+1}}) - \phi(X^\pi_{t_{i+1}})|^2\right]
	\wedge L'_f)
	\right] \\
	&
	+ \hat{\mathbb{E}}^x_{t_k, S}\left[\expectation{|f^{R}_{i+1}(X^\pi_{t_{i+1}}) - \tilde{f}^{R}_{i+1}(X^\pi_{t_i}, X^\pi_{t_{i+1}})|^2(1 - {\bf 1}_{\mathcal{A}}(X^\pi_{t_{i+1}}))}{k}{x}\right]
	=: \Xi^x_{t_k}(i,f),
\end{align*}
where $L'_y = 2LQ C^2_{p, A} + 2C_{Y,A}^2$ and $L'_f = 2LQ C^2_{p, A} + 2C^2_{f,A}$.
Note that although the size of $C_\pi:= \max_{k = 0, \ldots, N-1}\Delta_k$ may affect multiple constants here, like $C_{M, A}$ due to the probability law, $C_{Y, \mathcal{A}}$ and $C_{f, \mathcal{A}}$ by definition, however, $C_\pi\rightarrow 0$ would not make these constants converge to $0$.
So it may be easier to replace the constant $C_\pi$ by $T$ and consider these constants independent of the discretization scheme.
Therefore, we consider the refined regression error to be independent of the discretization scheme.

The following proposition summarizes the error bound for our scheme:
\begin{align*}
\Delta z^{}_k(x) := z^{}_k(x) - z^{R}_k(x);\quad
\Delta y^{}_k(x) := y^{}_k(x) - y^{R}_k(x).
\end{align*}
\begin{proposition}\label{proposition_a_s_upperbound} 
	Given Assumption \ref{assumption:globally_lipschitz}, and the time-grid $\pi$ and an $N$-dimensional vector $\gamma \in (0, +\infty)^N$ satisfying $12q(L_f^2R_\pi\vee 1)(\Delta_k + \frac{1}{\gamma_k})\leq 1$, for all $k\leq N-1$, we have, for $0 \leq k \leq N$,  
	\begin{align}
	& \hat{\mathbb{E}}^x_{t_k, S}[|\Delta y_k(x)|^2]\nonumber\\ 
	\leq &
	6q e^{T/4} \sum^{N-2}_{i=k} (\Delta_i + \gamma^{-1}_i) \Gamma_i L^2_f
	\Xi^x_{t_k}(i+1, y) + 3 e^{T/4}\sum^{N-1}_{i=k}(\Delta_i + \gamma^{-1}_i) \Gamma_i \frac{1}{\Delta_i} \Xi^x_{t_k}(i, y)\nonumber\\
	& + 3 e^{T/4} \sum^{N-1}_{i=k}(\Delta_i + \gamma^{-1}_i) \Gamma_i \Delta_i \Xi^x_{t_k}(i, f),\label{equation_y_error}
	\end{align}
	where $\Gamma_i := \prod^{k-1}_{i =0}(1+\gamma_i\Delta_i)$, and
	\begin{align*}
	& \hat{\mathbb{E}}^x_{t_k, S}\left[\sum^{N-1}_{i=k}\Delta_i\expectation{|\Delta z_i(X^\pi_{t_i})|^2}{k}{x}\Gamma_i \right]\\
	\leq & 
	(12q + 3Te^{T/4})\sum^{N-1}_{i = k+1} \left(\Delta_i + \gamma^{-1}_i\right) 
	\frac{1}{\Delta_i}\Xi^x_{t_k}(i, y)\Gamma_i
	+ 6q Te^{T/4} \sum^{N-2}_{i=k} (\Delta_i + \gamma^{-1}_i) \Gamma_i L^2_f \Xi^x_{t_k}(i+1, y)\\
	& + (12q + 3Te^{T/4})\sum^{N-1}_{i = k+1} \left(\Delta_i + \gamma^{-1}_i\right) \Delta_i 
	\Xi^x_{t_k}(i, f)\Gamma_i
	+ 4 \sum^{N-1}_{i=k} q \Xi^x_{t_k}(i, y)\Gamma_i.
	\end{align*}
\end{proposition} 
We will discuss the bound on $\Delta y_k$ here only as the two bounds are quite similar in structure. 
Note that the three terms within the sum at the right hand side of Equation \eqref{equation_y_error} are also of similar structure. 
They all sum up the refined regression error multiplied by some constant related to $\Delta_i$.
The most problematic term is $3 e^{T/4}\sum^{N-1}_{i=k}(\Delta_i + \gamma^{-1}_i) \Gamma_i \frac{1}{\Delta_i} \Xi^x_{t_k}(i, y)$ as the coefficient is $\mathcal{O}(1)$. 
The value $\sum^{N-1}_{i=k}(\Delta_i + \gamma^{-1}_i) \Gamma_i \frac{1}{\Delta_i}$ tends to infinity as the number of time steps tends to infinity. 
Therefore, one must use the parameters $M$ and $B$ to ensure the refined regression term is bounded by $C\Delta_i^{1+\epsilon}$ for some constant $C$, such that the sum and the error are bounded by $C C^\epsilon_\pi$.
This error plus the discretization error between the continuous system and the discretized system would be the complete error.
So, in practice, one should ensure that $N, M, M/B$ all tend together to infinity.

\begin{proof}
	The proof is fairly similar to the one used in \cite{gobet_turkedjiev_2016} with the necessary modifications for our present algorithm. 
	
	We shall derive an a-priori estimate of the error propagation in the recursion steps and we start with an estimate of $\Delta z^{}_k(x)$.
	Note that we add an extra term in the formula which is equal to zero due to the expectation of the Brownian motion being equal to zero.
	This term is added here to facilitate future steps of the proof.
	We have
	\begin{align*}
	|\Delta_k \Delta z^{}_k(x)|^2 = &  \left(\expectation{
		\left(
		\Delta y^{}_{k+1}(X^\pi_{t_{k+1}})-\expectation{\Delta y^{}_{k+1}(X^\pi_{t_{k+1}})}{k}{x}
		\right)
		\transpose{\Delta W_k}}{k}{x} \right.\\
	& + \left.
	\expectation{\left(y^{R}_{k+1}(X^\pi_{t_{k+1}})- \tilde{y}^{R}_{k+1}(X^\pi_{t_k}, X^\pi_{t_{k+1}})\right)\transpose{\Delta W_k}}{k}{x}\right)^2\\
	\leq & 2 \left(\expectation{
		\left(
		\Delta y^{}_{k+1}(X^\pi_{t_{k+1}})-\expectation{\Delta y^{}_{k+1}(X^\pi_{t_{k+1}})}{k}{x}
		\right)
		\transpose{\Delta W_k}}{k}{x} \right)^2\\
	& + 2 \left(\expectation{\left(y^{R}_{k+1}(X^\pi_{k+1})- \tilde{y}^{R}_{k+1}(X^\pi_{t_k}, X^\pi_{t_{k+1}})\right)\transpose{\Delta W_k}}{k}{x}\right)^2.
	\end{align*}
	The inequality follows from the inequality $(\sum^N_{n=1} a_n)^2 \leq \sum^N_{n=1} Na^2_n$, which will be frequently used in the proof and will not be specified again.
	By applying the Cauchy-Schwarz inequality, we can derive bounds for the two terms separately, where
	\begin{align*}
	& \left|
	\expectation{
		\left(
		\Delta y^{}_{k+1}(X^\pi_{t_{k+1}})-\expectation{\Delta y^{}_{k+1}(X^\pi_{t_{k+1}})}{k}{x}
		\right)
		\transpose{\Delta W_k}}{k}{x}
	\right|^2 \\
	\leq & 
	q \Delta_k 
	\left(
	\expectation{(\Delta y^{}_{k+1}(X^\pi_{t_{k+1}}))^2}{k}{x} - \left(\expectation{\Delta y^{}_{k+1}(X^\pi_{t_{k+1}})}{k}{x}\right)^2
	\right),
	\end{align*}
	and 
	\begin{align*}
	& \left|
	\expectation{
		\left(y^{R}_{k+1}(X^\pi_{k+1})- \tilde{y}^{R}_{k+1}(X^\pi_{t_k}, X^\pi_{t_{k+1}})\right)	
		\transpose{\Delta W_k}}{k}{x}
	\right|^2\\
	\leq & q \Delta_k 
	\expectation{
		\left|y^{R}_{k+1}(X^\pi_{k+1})- \tilde{y}^{R}_{k+1}(X^\pi_{t_k}, X^\pi_{t_{k+1}})\right|^2	
	}{k}{x}.
	\end{align*}
	Therefore, 
	\begin{align}
	\Delta_k |\Delta z^{}_k(x)|^2 \leq &
	2q \left(
	\expectation{(\Delta y^{}_{k+1}(X^\pi_{t_{k+1}}))^2}{k}{x} - \left(\expectation{\Delta y^{}_{k+1}(X^\pi_{t_{k+1}})}{k}{x}\right)^2
	\right)\nonumber\\
	& + 2q \expectation{
		\left|y^{R}_{k+1}(X^\pi_{k+1})- \tilde{y}^{R}_{k+1}(X^\pi_{t_k}, X^\pi_{t_{k+1}})\right|^2	
	}{k}{x}. \label{equation_estimates_for_z}
	\end{align}
	
	Combining the fact that $(a+b)^2 \leq (1 + \gamma_k\Delta_k)a^2 +(1 + \gamma_k^{-1}\Delta^{-1}_k) b^2$ for $(a,b)\in\mathbb{R}^2$, $\gamma_k >0$, and the Lipschitz property of $f$, one deduces with Equation (\ref{equation_estimates_for_z}) that, for $0\leq k\leq N-2$: 
	\begin{align}
	|\Delta y^{}_k(x)|^2 
	\leq &
	\left(
	\expectation{\Delta y^{}_{k+1}(X^\pi_{t_{k+1}})}{k}{x}
	+\expectation{y^{R}_{k+1}(X^\pi_{t_{k+1}}) - \tilde{y}^{R}_{k+1}(X^\pi_{t_k}, X^\pi_{t_{k+1}})}{k}{x}
	\right.\nonumber\\
	& + \expectation{f_{k+1}(y^{}_{k+1}(X^\pi_{t_{k+1}}), z^{}_{k+1}(X^\pi_{t_{k+1}})) - f^{R}_{k+1}(X^\pi_{t_{k+1}})}{k}{x}\Delta_k\nonumber\\
	& \left.
	+\expectation{
		f_{k+1}^{R}(X^\pi_{t_{k+1}})
		- \tilde{f}_{k+1}^{R}(X^\pi_{t_k}, X^\pi_{t_{k+1}})
	}{k}{x}\Delta_k
	\right)^2\nonumber\\
	\leq & 
	(1+ \gamma_k \Delta_k)\left(\expectation{\Delta y^{}_{k+1}(X^\pi_{t_{k+1}})}{k}{x}\right)^2 \nonumber\\
	& + 3\left(\Delta_k + \gamma^{-1}_k\right) \Delta_k 
	\left[
	L^2_f\expectation{(\Delta y^{}_{k+1}(X^\pi_{t_{k+1}}))^2}{k}{x} 
	+ L^2_f\expectation{(\Delta z^{}_{k+1}(X^\pi_{t_{k+1}}))^2}{k}{x}
	\right.\nonumber\\
	& \hspace{95pt}\left.
	+\frac{1}{\Delta_k^2}\expectation{|y^{}_{k+1}(X^\pi_{t_{k+1}}) - \tilde{y}^{}_{k+1}(X^\pi_{t_k}, X^\pi_{t_{k+1}})|^2}{k}{x}
	\right.\nonumber\\
	& \hspace{95pt}\left.
	+\expectation{
		|f_{k+1}^{R}(X^\pi_{t_{k+1}})
		- \tilde{f}_{k+1}^{R}(X^\pi_{t_k}, X^\pi_{t_{k+1}})|^2
	}{k}{x}
	\right]\nonumber\\
	\leq &
	(1+ \gamma_k \Delta_k)\left(\expectation{\Delta y^{}_{k+1}(X^\pi_{t_{k+1}})}{k}{x}\right)^2 \nonumber\\
	& + 3(\Delta_k + \gamma^{-1}_k) \Delta_k L^2_f\expectation{(\Delta y^{}_{k+1}(X^\pi_{t_{k+1}}))^2}{k}{x}\nonumber\\
	& + 6q(\Delta_k + \gamma^{-1}_k)  L^2_f R_\pi\left(
	\expectation{(\Delta y^{}_{k+2}(X^\pi_{t_{k+2}}))^2}{k}{x} - \expectation{\left(\expectation{\Delta y^{}_{k+2}(X^\pi_{t_{k+2}})}{k+1}{}\right)^2}{k}{x}
	\right) \nonumber\\
	& + 6q(\Delta_k + \gamma^{-1}_k)  L^2_f \expectation{
		\left|y^{R}_{k+2}(X^\pi_{t_{k+2}})- \tilde{y}^{R}_{k+2}(X^\pi_{t_{k+1}}, X^\pi_{t_{k+2}})\right|^2	
	}{k}{x}\nonumber\\
	& + 3(\Delta_k + \gamma^{-1}_k) \Delta_k \frac{1}{\Delta_k^2}\expectation{|y^{R}_{k+1}(X^\pi_{t_{k+1}}) - \tilde{y}^{R}_{k+1}(X^\pi_{t_k}, X^\pi_{t_{k+1}})|^2}{k}{x}\nonumber\\
	& + 3(\Delta_k + \gamma^{-1}_k) \Delta_k \expectation{
		|f_{k+1}^{R}(X^\pi_{t_{k+1}})
		- \tilde{f}_{k+1}^{R}(X^\pi_{t_k}, X^\pi_{t_{k+1}})|^2
	}{k}{x}, \label{equation_estimation_y_lipschitz}
	\end{align}
	while 
	\begin{align}
	|\Delta y^{}_{N-1}(x)|^2 
	\leq & 
	3\left(\Delta_k + \gamma^{-1}_k\right) \Delta_k 
	\left[
	\frac{1}{\Delta_k^2}\expectation{|y^{R}_{k+1}(X^\pi_{t_{k+1}}) - \tilde{y}^{R}_{k+1}(X^\pi_{t_k}, X^\pi_{t_{k+1}})|^2}{k}{x}
	\right.\nonumber\\
	& \hspace{95pt}\left.
	+\expectation{
		|f_{k+1}^{R}(X^\pi_{t_{k+1}})
		- \tilde{f}_{k+1}^{R}(X^\pi_{t_k}, X^\pi_{t_{k+1}})|^2
	}{k}{x}
	\right]. \label{equation_estimation_for_y}
	\end{align}
	Next, we define the following sequence 
	\begin{equation*}
	\lambda_k := \left[1 + \left(\gamma_{k-1} + \frac{1}{4}\right)\Delta_{k-1}\right]\lambda_{k-1}, \mbox{ where } \lambda_0:=1,
	\end{equation*}
	consider the sum of $|\Delta y_i(X^\pi_{t_i})|\lambda_i$, from $i = 1$ to $N-1$, and take conditional expectations with respect to $\mathcal{F}_k$.
	Applying Equation (\ref{equation_estimation_for_y}) for the case $k = N-1$ and Equation (\ref{equation_estimation_y_lipschitz}) otherwise, we have:
	\begin{align*}
	\sum^{N-1}_{i=k} \expectation{|\Delta y^{}_i(X^\pi_{t_{i}})|^2 \lambda_i}{k}{x}
	\leq & \sum^{N-2}_{i=k}\lambda_{i+1}\expectation{\left(\Delta y^{}_{i+1}(X^\pi_{t_{i+1}})\right)^2}{k}{x}\\
	& + \sum^{N-2}_{i=k} 6q(\Delta_i + \gamma^{-1}_i)  L^2_f \lambda_i \expectation{
		\left|y^{R}_{i+2}(X^\pi_{i+2})- \tilde{y}^{R}_{i+2}(X^\pi_{i+1}, X^\pi_{i+2})\right|^2	
	}{k}{x}\\
	& + \sum^{N-1}_{i=k}3(\Delta_i + \gamma^{-1}_i) \Delta_i \frac{1}{\Delta_i^2}\lambda_i \expectation{|y^{R}_{i+1}(X^\pi_{t_{i+1}}) - \tilde{y}^{R}_{i+1}(X^\pi_{t_i}, X^\pi_{t_{i+1}})|^2}{k}{x}\\
	& + \sum^{N-1}_{i=k}3(\Delta_i + \gamma^{-1}_i) \Delta_i \lambda_i\expectation{
		|f_{i+1}^{R}(X^\pi_{t_{i+1}})
		- \tilde{f}_{i+1}^{R}(X^\pi_{t_{i}}, X^\pi_{t_{i+1}})|^2
	}{k}{x}.
	\end{align*}
	By rearranging the terms, we have:
	\begin{align*}
	|\Delta y^{}_k(x)|^2\lambda_k
	\leq &
	\sum^{N-2}_{i=k} 6q(\Delta_i + \gamma^{-1}_i)  L^2_f \lambda_i \expectation{
		\left|y^{R}_{i+2}(X^\pi_{t_{i+2}})- \tilde{y}^{R}_{i+2}(X^\pi_{t_{i+1}}, X^\pi_{t_{i+2}})\right|^2	
	}{k}{x}\\
	& + \sum^{N-1}_{i=k}3(\Delta_i + \gamma^{-1}_i) \Delta_i \frac{1}{\Delta_i^2}\lambda_i \expectation{|y^{R}_{i+1}(X^\pi_{t_{i+1}}) - \tilde{y}^{R}_{i+1}(X^\pi_{t_i}, X^\pi_{t_{i+1}})|^2}{k}{x}\\
	& + \sum^{N-1}_{i=k}3(\Delta_i + \gamma^{-1}_i) \Delta_i \lambda_i\expectation{
		|f_{i+1}^{R}(X^\pi_{t_{i+1}})
		- \tilde{f}_{i+1}^{R}(X^\pi_{t_{i}}, X^\pi_{t_{i+1}})|^2
	}{k}{x}.
	\end{align*}
	It follows from the simple inequality $\Gamma_k \leq \lambda_k = \exp(\sum^k_{i=0}\log(1+(\gamma_i +0.25)\Delta_i) \leq e^{T/4}\Gamma_k$ that, for all $k\in \{0, \ldots, N\}$,
	\begin{align}
	|\Delta y^{}_k(x)|^2 
	\leq & 
	|\Delta y^{}_k(x)|^2\Gamma_k \nonumber \\
	\leq &
	6q e^{T/4} \sum^{N-2}_{i=k} (\Delta_i + \gamma^{-1}_i) \Gamma_i L^2_f \expectation{
		\left|y^{R}_{i+2}(X^\pi_{t_{i+2}})- \tilde{y}^{R}_{i+2}(X^\pi_{t_{i+1}}, X^\pi_{t_{i+2}})\right|^2	
	}{k}{x}\nonumber\\
	&
	+ 3 e^{T/4}\sum^{N-1}_{i=k}(\Delta_i + \gamma^{-1}_i) \Gamma_i \frac{1}{\Delta_i} \expectation{|y^{R}_{i+1}(X^\pi_{t_{i+1}}) - \tilde{y}^{R}_{i+1}(X^\pi_{t_i}, X^\pi_{t_{i+1}})|^2}{k}{x}\nonumber\\
	& 
	+ 3 e^{T/4} \sum^{N-1}_{i=k}(\Delta_i + \gamma^{-1}_i) \Gamma_i \Delta_i \expectation{
		|f_{i+1}^{R}(X^\pi_{t_{i+1}})
		- \tilde{f}_{i+1}^{R}(X^\pi_{t_{i}}, X^\pi_{t_{i+1}})|^2
	}{k}{x}. \label{equation_pointwise_bound_y}
	\end{align}
	
	We can take expectations with respect to the simulation cloud and apply Theorem \ref{theorem_regression}, which finishes the calculation for $\Delta y$.
	
	Regarding the error term $\Delta z$,  $\sum^{N-1}_{i=k}\Delta_i\expectation{|\Delta z^{}_i(X^\pi_{t_i})|^2}{k}{x}\Gamma_i $ is bounded from above by
	\begin{align*}
	& \sum^{N-1}_{i=k}\Delta_i\expectation{|\Delta z^{}_i(X^\pi_{t_i})|^2}{k}{x}\Gamma_i \\
	\leq &
	\sum^{N-1}_{i=k} 2q \left(
	\expectation{(\Delta y^{}_{i+1}(X^\pi_{t_{i+1}}))^2}{k}{x} - \expectation{\left(\expectation{\Delta y^{}_{i+1}(X^\pi_{t_{i+1}})}{i}{}\right)^2}{k}{x} 
	\right)\Gamma_{i+1}\\
	& + \sum^{N-1}_{i=k} 2q \expectation{
		\left|y^{R}_{i+1}(X^\pi_{t_{i+1}})- \tilde{y}^{R}_{i+1}(X^\pi_{t_i}, X^\pi_{t_{i+1}})\right|^2	
	}{k}{x}\Gamma_i\\
	\leq & 
	2q \Gamma_N \expectation{(\Delta y^{}_N(X^\pi_{t_N}))^2}{k}{x}\\
	& + \sum^{N-1}_{i=k+1} 2q \Gamma_i \left(
	\expectation{(\Delta y^{}_i(X^\pi_{t_i}))^2}{k}{x} - (1+\gamma_i \Delta_i)\expectation{\left(\expectation{\Delta y^{}_{i+1}(X^\pi_{t_{i+1}})}{i}{}\right)^2}{k}{x} 
	\right)\\
	& + \sum^{N-1}_{i=k} 2q \expectation{
		\left|y^{R}_{i+1}(X^\pi_{t_{i+1}})- \tilde{y}^{R}_{i+1}(X^\pi_{t_i}, X^\pi_{t_{i+1}})\right|^2	
	}{k}{x}\Gamma_i,
	\end{align*}
	because of Equation (\ref{equation_estimates_for_z}), and from (\ref{equation_estimation_y_lipschitz}), we have
	\begin{align*}
	\sum^{N-1}_{i=k}\Delta_i\expectation{|\Delta z^{}_i(X^\pi_{t_i})|^2}{k}{x}\Gamma_i
	\leq & 
	6\sum^{N-1}_{i = k+1} q \left(\Delta_i + \gamma^{-1}_i\right) \Delta_i 
	L^2_f\expectation{(\Delta y^{}_{i+1}(X^\pi_{t_{i+1}}))^2}{k}{x} \Gamma_i \\
	& + 6\sum^{N-1}_{i = k+1} q \left(\Delta_i + \gamma^{-1}_i\right) \Delta_i 
	L^2_f\expectation{(\Delta z^{}_{i+1}(X^\pi_{t_{i+1}}))^2}{k}{x} \Gamma_i\\
	& + 6\sum^{N-1}_{i = k+1} q \left(\Delta_i + \gamma^{-1}_i\right) 
	\frac{1}{\Delta_i}\expectation{|y^{R}_{i+1}(X^\pi_{t_{i+1}}) - \tilde{y}^{R}_{i+1}(X^\pi_{t_i}, X^\pi_{t_{i+1}})|^2}{k}{x}\Gamma_i\\
	& + 6\sum^{N-1}_{i = k+1} q \left(\Delta_i + \gamma^{-1}_i\right) \Delta_i 
	\expectation{
		|f_{i+1}^{R}(X^\pi_{t_{i+1}})
		- \tilde{f}_{i+1}^{R}(X^\pi_{t_{i}}, X^\pi_{t_{i+1}})|^2
	}{k}{x}\Gamma_i\\
	& + \sum^{N-1}_{i=k} 2q \expectation{
		\left|y^{R}_{i+1}(X^\pi_{i+1})- \tilde{y}^{R}_{i+1}(X^\pi_{t_i}, X^\pi_{t_{i+1}})\right|^2	
	}{k}{x}\Gamma_i.
	\end{align*}
	
	Using the assumptions of the proposition statement, it follows that
	\begin{align*}
	& \sum^{N-1}_{i=k}\Delta_i\expectation{|\Delta z^{}_i(X^\pi_{t_i})|^2}{k}{x}\Gamma_i \\
	\leq & 
	12\sum^{N-1}_{i = k+1} q \left(\Delta_i + \gamma^{-1}_i\right) 
	\frac{1}{\Delta_i}\expectation{|y^{R}_{i+1}(X^\pi_{t_{i+1}}) - \tilde{y}^{R}_{i+1}(X^\pi_{t_i}, X^\pi_{t_{i+1}})|^2}{k}{x}\Gamma_i\\
	& + 12\sum^{N-1}_{i = k+1} q \left(\Delta_i + \gamma^{-1}_i\right) \Delta_i 
	\expectation{
		|f_{i+1}^{R}(X^\pi_{t_{i+1}})
		- \tilde{f}_{i+1}^{R}(X^\pi_{t_{i}}, X^\pi_{t_{i+1}})|^2
	}{k}{x}\Gamma_i\\
	& + 4 \sum^{N-1}_{i=k} q \expectation{
		\left|y^{R}_{i+1}(X^\pi_{i+1})- \tilde{y}^{R}_{i+1}(X^\pi_{t_i}, X^\pi_{i+1})\right|^2	
	}{k}{x}\Gamma_i
	+ \sum^{N-1}_{j = k+1} \Delta_j \expectation{(\Delta y^{}_{j+1}(X^\pi_{t_{j+1}}))^2}{k}{x} \Gamma_{j+1}.
	\end{align*}
	Note that we may bound each individual term in the last sum with the estimate from Equation (\ref{equation_pointwise_bound_y}) and by taking conditional expectations.
	\begin{align*}
	\mathbb{E}^x_{t_k}[(\Delta y_{j+1}(X^\pi_{t_{j+1}}))^2 \Gamma_{i+1}]
	\leq &
	6q e^{T/4} \sum^{N-2}_{i=j+1} (\Delta_i + \gamma^{-1}_i) \Gamma_i L^2_f \expectation{
		\left|y^{R}_{i+2}(X^\pi_{t_{i+2}})- \tilde{y}^{R}_{i+2}(X^\pi_{t_{i+1}}, X^\pi_{t_{i+2}})\right|^2	
	}{k}{x}\nonumber\\
	& + 3 e^{T/4}\sum^{N-1}_{i=j+1}(\Delta_i + \gamma^{-1}_i) \Gamma_i \frac{1}{\Delta_i} \expectation{|y^{R}_{i+1}(X^\pi_{t_{i+1}}) - \tilde{y}^{R}_{i+1}(X^\pi_{t_i}, X^\pi_{t_{i+1}})|^2}{k}{x}\nonumber\\
	& + 3 e^{T/4} \sum^{N-1}_{i=j+1}(\Delta_i + \gamma^{-1}_i) \Gamma_i \Delta_i \expectation{
		|f_{i+1}^{R}(X^\pi_{t_{i+1}})
		- \tilde{f}_{i+1}^{R}(X^\pi_{t_{i}}, X^\pi_{t_{i+1}})|^2
	}{k}{x}\\
	\leq &
	6q e^{T/4} \sum^{N-2}_{i=k+1} (\Delta_i + \gamma^{-1}_i) \Gamma_i L^2_f \expectation{
		\left|y^{R}_{i+2}(X^\pi_{t_{i+2}})- \tilde{y}^{R}_{i+2}(X^\pi_{t_{i+1}}, X^\pi_{t_{i+2}})\right|^2	
	}{k}{x}\nonumber\\
	& + 3 e^{T/4}\sum^{N-1}_{i=k+1}(\Delta_i + \gamma^{-1}_i) \Gamma_i \frac{1}{\Delta_i} \expectation{|y^{R}_{i+1}(X^\pi_{t_{i+1}}) - \tilde{y}^{R}_{i+1}(X^\pi_{t_i}, X^\pi_{t_{i+1}})|^2}{k}{x}\nonumber\\
	& + 3 e^{T/4} \sum^{N-1}_{i=k+1}(\Delta_i + \gamma^{-1}_i) \Gamma_i \Delta_i \expectation{
		|f_{i+1}^{R}(X^\pi_{t_{i+1}})
		- \tilde{f}_{i+1}^{R}(X^\pi_{t_{i}}, X^\pi_{t_{i+1}})|^2
	}{k}{x}
	\end{align*}
	This upper bound is independent of $j$.
	Summing up the remaining parts, the time increments, results in the full time length $T$.
	We have:
	\begin{align*}
	& \sum^{N-1}_{i=k}\Delta_i\expectation{|\Delta z^{}_i(X^\pi_{t_i})|^2}{k}{x}\Gamma_i \\
	\leq & 
	(12q + 3Te^{T/4})\sum^{N-1}_{i = k+1} \left(\Delta_i + \gamma^{-1}_i\right) 
	\frac{1}{\Delta_i}\expectation{|y^{R}_{i+1}(X^\pi_{t_{i+1}}) - \tilde{y}^{R}_{i+1}(X^\pi_{t_i}, X^\pi_{t_{i+1}})|^2}{k}{x}\Gamma_i\\
	& + 6q Te^{T/4} \sum^{N-2}_{i=k} (\Delta_i + \gamma^{-1}_i) \Gamma_i L^2_f \expectation{
		\left|y^{R}_{i+2}(X^\pi_{t_{i+2}})- \tilde{y}^{R}_{i+2}(X^\pi_{t_{i+1}}, X^\pi_{i+2})\right|^2	
	}{k}{x}\\
	& + (12q + 3Te^{T/4})
	\sum^{N-1}_{i = k+1} \left(\Delta_i + \gamma^{-1}_i\right) \Delta_i 
	\expectation{
		|f_{i+1}^{R}(X^\pi_{t_{i+1}})
		- \tilde{f}_{i+1}^{R}(X^\pi_{t_{i}}, X^\pi_{t_{i+1}})|^2
	}{k}{x}\Gamma_i\\
	& + 4 \sum^{N-1}_{i=k} q \expectation{
		\left|y^{R}_{i+1}(X^\pi_{i+1})- \tilde{y}^{R}_{i+1}(X^\pi_{t_i}, X^\pi_{t_{i+1}})\right|^2	
	}{k}{x}\Gamma_i,
	\end{align*}
	Again, taking expectations with respect to the simulation cloud finishes the proof.
\end{proof}

\section{Numerical Experiments} \label{section_numerical}

In this section, numerical experiments are conducted for some selected examples.
Before discussing these examples, we would specify the forward and backward discretization scheme used in these experiments.
In particular, we introduce a more general backward scheme to show that our algorithm can be applied in general circumstances.

\subsection{Forward and Backward Scheme}
In this section, we conduct our numerical experiments with the Euler-Maruyama discretization scheme, which is a common standard in the literature.

\begin{definition}[Euler-Maruyama scheme]
	The Euler-Maruyama scheme is defined by 
	\begin{equation*}
	X^\pi_{t_{k+1}} = X^\pi_{t_k} + b(t_k,X^\Delta_{t_k})\Delta_k + \sigma(t_k,X^\pi_{t_k})\Delta W_k =: d(X^\pi_{t_k}, \Delta W_k).
	\end{equation*}
	
	Note that the conditional expectation $\expectation{\frac{\Delta W_{l,k}}{\Delta_k}p(X^\pi_{t_{k+1}})}{k}{x}$ can be calculated by:
	\begin{align*}
	\expectation{\frac{\Delta W_{l,k}}{\Delta_k}p(X^\pi_{t_{k+1}})}{k}{x}
	= & \frac{1}{\sqrt{(2\pi)^q \Delta^q_k}} 
	\int_{\mathbb{R}^q}p(d(x,y))
	\frac{\partial}{\partial y_l}\left(-\exp\left(-\frac{1}{2}\sum^q_{r=1}\frac{y^2_r}{\Delta_k}\right)\right)dy\\
	= & \frac{1}{\sqrt{(2\pi)^q\Delta^q_k}}\int_{\mathbb{R}^q}\exp\left(-\frac{1}{2}\sum^q_{r=1}\frac{y^2_r}{\Delta_p}\right) \nabla p(d(x,y))\frac{\partial d(x,y)}{\partial y_l}dy\\
	= & \expectation{\nabla p(X^\pi_{t_{k+1}})}{k}{x}\sigma_{l}(t_k,x),
	\end{align*}
	where $\sigma_l$ is the $l$-th column of the matrix $\sigma$. 
\end{definition}

For example, for the one-dimensional monomial $x^r$, $r \in \mathbb{N}$ and a forward process discretized by the Euler-Maruyama scheme, we have
\begin{align*}
\expectation{\frac{\Delta W_k}{\Delta_k}(X^\pi_{t_{k+1}})^r}{k}{x}
= & \expectation{r(X^\pi_{t_{k+1}})^{r-1}}{k}{x}\sigma(t_k,x).
\end{align*}
The conditional expectations of polynomials are calculated directly by definition. 
We have 
\begin{align*}
\mathbb{E}^x_{t_k}[(X^\pi_{t_{k+1}})^0] = & 1;\\
\mathbb{E}^x_{t_k}[(X^\pi_{t_{k+1}})^1] = & x + b(t_k,x)\Delta_k;\\
\mathbb{E}^x_{t_k}[(X^\pi_{t_{k+1}})^{2}] = & x^2 + 2 x b(t_k,x)\Delta_k + \sigma(t_k, x)^2\Delta_k + b(t_k, x)^2 \Delta^2_k ,\\
\end{align*}
and so on.

For the backward discretization, we apply the theta-scheme from \cite{zhao_li_zhang_2012} and \cite{ruijter_oosterlee_2015}:
\begin{align*}
Y^{\pi}_{t_N} = & \Phi(X^\pi_{t_N}),\quad Z^{\pi}_{t_N} = \transpose{\left(\nabla \Phi(X^\pi_{t_N}) \sigma (t_N, X^\pi_{t_N})\right)},\\
Z^{\pi}_{t_k} = & -\theta_2^{-1} (1-\theta_2) \expectation{Z^\pi_{t_{k+1}}}{k}{} + \frac{1}{\Delta_k}\theta_2^{-1}\expectation{Y^\pi_{t_{k+1}}\Delta W_k}{k}{} \\
& + \theta_2^{-1}(1-\theta_2)\expectation{f_{k+1}(Y^\pi_{t_{k+1}}, Z^\pi_{t_{k+1}})\Delta W_k }{k}{}, \;
k = N-1, \ldots, 0,\\
Y^\pi_{t_k} = & \expectation{Y^\pi_{t_{k+1}}}{k}{} + \Delta_k \theta_1 f_k(Y^\pi_{t_k},  Z^\pi_{t_k})\\
& + \Delta_k (1-\theta_1) \expectation{f_{k+1}(Y^\pi_{t_{k+1}},Z^\pi_{t_{k+1}})}{k}{},\; k= N-1,\ldots,0,
\end{align*}
$0 \leq \theta_1 \leq 1$ and $0 < \theta_2 \leq 1$.

By picking various parameters $(\theta_1, \theta_2)$, we can construct different types of one-step dynamic programming schemes.
For example, the choice $(\theta_1, \theta_2) = (0,1)$ would result in an explicit scheme, while the choice $(\theta_1, \theta_2) = (0.5, 0.5)$ would give the Crank-Nicolson scheme.
Using a general construction means that our algorithm can be applied to various types of schemes and we may adjust our algorithms towards the specific problem.

Applying the SGBM algorithm to this general scheme, we have that the approximate functions within the bundle at time $k$ are defined by:
\begin{align*}
& \approximant{z}{r, k}{(\theta_1, \theta_2), R}(b, x) = 
-\theta^{-1}_2(1-\theta_2)\expectation{p(X^\pi_{t_{k+1}})}{k}{x}
\regressionparameter{\beta}{r, k+1}\\
& \hspace{80pt} 
+ \theta_2^{-1}\expectation{
	\frac{\Delta W_{r,k}}{\Delta_k}p(X^\pi_{t_{k+1}})
}{k}{x} 
(\regressionparameter{\alpha}{k+1}+(1-\theta_2)\Delta_k\regressionparameter{\gamma}{k+1}),\quad
r= 1,\ldots,d;\\
& \approximant{y}{k}{(\theta_1, \theta_2), R, 0}(b, x) = 
\expectation{p(X^\pi_{t_{k+1}})}{k}{x} \regressionparameter{\alpha}{k+1},\\
& \approximant{y}{k}{(\theta_1, \theta_2), R, i}(b, x) =
\Delta_k\theta_1 f(t_k, x, \approximant{y}{k}{(\theta_1, \theta_2), R,i-1}(x), \approximant{z}{k}{(\theta_{1}, \theta_2), R}(x)) + \approximant{h}{k}{}(b, x),\\
& \approximant{h}{k}{}(b, x) = 
\expectation{p(X^\pi_{t_{k+1}})}{k}{x}
(\regressionparameter{\alpha}{k+1}+ \Delta_k(1 - \theta_1)\regressionparameter{\gamma}{k+1}),
\quad i = 1, \ldots, I.
\end{align*}
Note that a Picard iteration is performed at each time step for each bundle if the choice of $(\theta_1, \theta_2)$ results in an implicit scheme. 
For further details on the application of the Picard iteration, readers may refer to \cite{gobet_lemor_warin_2005} or \cite{ruijter_oosterlee_2015} and the references therein.

Different types of backward discretizations will be considered for Example 1.

\subsection{Example 1}
This example is originally from \cite{zhao_li_zhang_2012}. The considered FBSDE is given by
\begin{equation*}
\left\{
\begin{array}{l}
dX_t =  d\omega_t,\\
dY_t = -(Y_tZ_t-Z_t+2.5Y_t-\sin(t+X_t)\cos(t+X_t)-2\sin(t+X_t))dt+Z_td\omega_t.
\end{array}
\right.
\end{equation*}
We take the initial and terminal conditions $x_0 = 0$ and $Y_T = \sin(X_T+T)$.

The exact solution is given by 
\begin{equation*}
(Y_t,Z_t)=(\sin(X_t+t),\cos(X_t+t)).
\end{equation*}
The terminal time is set to be $T=1$ and $(Y_0,Z_0)=(0,1)$. 
We use the set $\{1, x, x^2\}$ as the regression base for this example. 
We apply equal partitioning bundling for all our tests with the sample paths sorted by the value function $x$.
As mentioned in Session \ref{session_assumptions}, not all assumptions set in this work are necessary for the basic SGBM algorithm to work. 
For example, Assumption \ref{assumption:globally_lipschitz} is included to ensure the existence and uniqueness of the solution of the BSDE.
In this example, even though the driver is not Lipschitz, one can check that the above solution solves the BSDE with It\^o's formula, and the SGBM algorithm still applies.

Table \ref{test_case_general} shows the tests that we have run.
Basically, our test cases can be placed into two groups.
Test cases 1a, 1b, 1c are tests for the explicit version of our algorithm, while test cases 1d, 1e, 1f are for the Crank-Nicolson version.
Within each group, the three tests are run for identical test settings, except for the constant $L$, i.e., the pre-set limit for the Euclidean norm so that we may check the influence of the factor $L$.
Within each test, the factors $M$, $N$ and $B$ are linked to a common factor $J$ such that when $J$ tends infinity, $N$, $B$ and $M/B$ tend to infinity as well.
This setting is inspired by our observation on the error bound that all three factors should tend to infinity together to ensure the convergence of the algorithm. 
However, the extra ratio between the three factors is from empirical experience.

\begin{table}[h]
	\centering
	\begin{tabular}{|c|ccccccc|}
		\hline
		Test Case & $\theta_1$ & $\theta_2$ & I & M & N & B & L\\
		\hline
		1a &  0 & 1 & - & $2^{2J}$ & $2^J$ & $2^J$ & $100$\\
		1b &  0 & 1 & - & $2^{2J}$ & $2^J$ & $2^J$ & $10000$\\
		1c &  0 & 1 & - & $2^{2J}$ & $2^J$ & $2^J$ & $-$\\
		1d &  0.5 & 0.5 & 4 & $2^{2J}$ & $2^J$ & $2^J$ & $100$\\
		1e &  0.5 & 0.5 & 4 & $2^{2J}$ & $2^J$ & $2^J$ & $10000$\\
		1f &  0.5 & 0.5 & 4 & $2^{2J}$ & $2^J$ & $2^J$ & $-$\\
		\hline
	\end{tabular}
	\caption{Test cases for Example 1} \label{test_case_general}
\end{table}

\subsection{Example 2: Black-Scholes European option}
The second example under consideration is the calculation of the price $v(t,S_t)$ of a European option under the  $d$-dimensional Black-Scholes model by solving a FBSDE, which has been a classical application of BSDEs.
It has been introduced in classical papers, like \cite{karoui1997backward}, and here we will provide a brief review.
We consider a market where the assets satisfy:
\begin{equation*}
dS_{i,t} = \mu_i S_{i,t} dt + \sigma_i S_{i,t} d B_{i, t},\; 1 \leq i \leq d,
\end{equation*}
where $B_t$ is a correlated $d$-dimensional Wiener process, with $$dB_{i, t}dB_{j, t} = \rho_ij dt.$$
The parameters $\rho_{ij}$ form a symmetric non-negative matrix $\rho$, 
\begin{equation*}
\rho = 
\left(\begin{array}{ccccc}
1 & \rho_{12} & \rho_{13} & \cdots & \rho_{1q}\\
\rho_{21} & 1 & \rho_{23} & \cdots & \rho_{2q}\\
\vdots & \vdots & \vdots & & \vdots\\
\rho_{q1} & \rho_{q2}& \rho_{q3} & \cdots & 1
\end{array}\right),
\end{equation*}
and we assume it is invertible.
By performing a Cholesky decomposition on $\rho$ such that $\mathfrak{C} \mathfrak{C}^\top = \rho$, where $\mathfrak{C}$ is a lower triangular matrix with real and positive diagonal entries, we may relate the correlated and standard Brownian motions, as follows, $$B_t = \mathfrak{C} W_t.$$

Along the line of reasoning in \cite{ruijter_oosterlee_2015}, we assume the financial market is complete, there is no trading restriction and a derivative can be perfectly hedged. 
To derive the corresponding pricing BSDE for a European option with terminal payoff $g(S_T)$, we construct a replicating portfolio $Y_t$, containing $\omega_{i, t}$ of asset $S_{i,t}$ and bonds with risk-free return rate $r$.
Applying the self-financing assumption, the portfolio follows the SDE:
\begin{equation*}
dY_t = -(-r Y_t - \sum^d_{i=1} \omega_{i,t}(\mu_i-r)S_{i,t})dt + \sum^d_{i=1}\omega_{i. t}\sigma_i S_i d B_{i, t}.
\end{equation*}
If we set $Z_t = (\omega_{1,t} \sigma_1 S_{1,t}, \ldots, \omega_{d, t}\sigma_d, S_{d, t})\mathfrak{C}$, then $(Y, Z)$ solves the BSDE,
\begin{equation*}
\begin{cases}
dY_t = -\left(-r Y_t - Z_t \mathfrak{C}^{-1}\left(\frac{\mu - r}{\sigma}\right) \right)dt + Z_t dW_t;\\
Y_T = g(S_T),
\end{cases}
\end{equation*}
where $\left(\frac{\mu -r}{\sigma}\right) = \left(\frac{\mu_1-r}{\sigma_1}, \cdots, \frac{\mu_q - r}{\sigma_q}\right)^T$.

We test our algorithm for the next two cases.

\subsubsection{Arithmetic Basket Put Option}

In this numerical test, we use the 5-dimensional example from \cite{2007reisingerefficient}, which is designed as a tractable representation for the German stock index DAX at that time.
All $\mu_i$ are assumed to be $r$ here.
The volatilities are given by $$(\sigma_1, \sigma_2, \sigma_3, \sigma_4, \sigma_5) = (0.518, 0.648, 0.623, 0.570, 0.530), $$ while the correlations $\rho$ are given by  
\begin{equation*}
\rho = 
\left(\begin{array}{ccccc}
1.00 & 0.79 & 0.82 & 0.91 & 0.84\\
0.79 & 1.00 & 0.73 & 0.80 & 0.76\\
0.82 & 0.73 & 1.00 & 0.77 & 0.72\\
0.91 & 0.80 & 0.77 & 1.00 & 0.90\\
0.84 & 0.76 & 0.72 & 0.90 & 1.00
\end{array}\right).
\end{equation*}

We would consider a European weighted basket put option for $T=1$ year, with the payoff function $g$ given by $$g(S) = \left(1 - \sum^5_{i=1}w_i S_i\right)^+,$$
where $(w_1, w_2, w_3, w_4, w_5) = (38.1, 6.5, 5.7, 27.0, 22.7)$.
The risk free interest rate is $r=0.05$ and all the stocks have starting value 0.01.
The reference price is given as 0.175866 in \cite{2007reisingerefficient}.

We perform the equal-partitioning bundling technique and sort the paths in different bundles according to the ordering of the value $\sum^5_{i = 1}w_i S^m_{i, t_p}$.
The regression basis is chosen to be $p_k(x) = \left(\sum^5_{i = 1}w_i x_i\right)^{k-1}$ for $k = 1, \ldots, K$.

Table \ref{test_case_arithmetic} shows the tests that we have run.
In these tests, we keep most of the parameters fixed but vary the number of bundles.
We test our algorithm for the explicit scheme with a second-order regression basis and the Crank-Nicolson scheme with a third-order regression basis.
The change of basis is made to test the impact of the regression basis to our algorithm.
We just keep these two sets of tests to demonstrate the impact of the number of bundles.

\begin{table}[h]
	\centering
	\begin{tabular}{|c|cccccccc|}
		\hline
		Test Case & $\theta_1$ & $\theta_2$ & I & M & N & B & L & K\\
		\hline
		2.1a & 0.5 & 0.5 & 4 & $2^{12}$ & $10$ & $2^{2J}$ & -& 3\\
		2.1b & 0 & 1 & - & $2^{12}$ & $10$ & $2^{2J}$ & -& 2\\
		\hline
	\end{tabular}
	\caption{Test cases for Example 2.1} \label{test_case_arithmetic}
\end{table}

\subsubsection{Example 2.2: Geometric Basket Put Option} 

Here we also consider the problem of pricing q-dimensional geometric basket options with initial state $S_0 = (40, \ldots, 40) \in \mathbb{R}^q$; strike $K = 40$; risk-free interest rate $r = 0.06$; volatility $\sigma_i = 0.2, i = 1, \ldots, d$; correlation $\rho_{ij} = 0.25, i,j = 1, \ldots, d, i\neq j$; and maturity $T=1.0$.
The final payoff function is given by $$g(S) = \left(K - \left(\prod^d_{i = 1} S_i\right)^{\frac{1}{d}}\right)^+.$$
This is the same setting as in \cite{leitao_oosterlee_2015} but for European options instead of Bermudan options.

We again use the equal-partitioning technique and sort the paths in different bundles according to the ordering of the values $\left(\prod^d_{i=1}S^m_{i, t_p}\right)^{\frac{1}{d}}$. 
The regression basis is chosen to be $p_k(x) = \left(\prod^d_{i = 1} x_i\right)^{\frac{k-1}{d}}$ for $k = 1, \ldots, 3$.

Since the geometric product of a geometric Brownian motion remains a geometric Brownian motion, the analytic solution can be found using Black-Scholes formula and any other classical pricing method.

Table \ref{test_case_geometric} shows the tests that we have run.
In these sets of tests, we fixed all the parameters but change the number of stocks in our test.
This example is used to test the {\em scalability} of our methodology.
Tests are performed for both explicit and Crank-Nicolson schemes.

\begin{table}[h]
	\centering
	\begin{tabular}{|c|ccccccc|}
		\hline
		Test Case & $\theta_1$ & $\theta_2$ & I & M & N & B & L\\
		\hline
		2.2a & 0 & 1 & - & $2^{12}$ & 20 & 16 & - \\
		2.2b & 0.5 & 0.5 & 4 & $2^{12}$ & 20 & 16 & - \\
		\hline
	\end{tabular}
	\caption{Test cases for Example 2.2} \label{test_case_geometric}
\end{table}

\subsection{Results}
The results are given as the average values of 10 separated runs of the algorithm.

We first consider the results of the explicit version of our algorithm applied to Example 1, namely test cases 1a, 1b and 1c, in Table \ref{test_result_explicit}.
This test can be seen as a proof of concept.
As mentioned, we design the test in such a way that the number of steps $N$, the number of bundles $B$  and the ratio $M/B$ all tend to infinity.
As expected, our algorithm converges under this setting.
Moreover, the total variation of the absolute errors among each successful run converges with respect to $J$ too, as the reader can read from the second part of Table \ref{test_result_explicit}.
It is defined as the sum of the individual differences between the Monte Carlo result of each run (which is not rejected) and the analytic solution, divided by the total number of successful runs.

\begin{table}[h]
	\centering
	\begin{tabular}{c|ccccccc}
		\multicolumn{8}{c}{$|Y_0-y^{(\theta_1, \theta_2), R}_0(x_0)|$}\\
		J & 2 & 3 & 4 & 5  & 6 & 7 & 8\\
		\hline
		1a & 0.023535 & 0.20392 & 0.046947 & 0.057056 & 0.026622 & 0.018172 & 0.016179\\
		1b & 0.18360 & 0.17807 & 0.098821 & 0.030159 & 0.028840 & 0.019621 & 0.0057568\\
		1c & 0.41648 & 0.14362 & 0.10368 & 0.04658 & 0.018068 & 0.019175 & 0.0098448\\
		\hline
		\multicolumn{8}{c}{Total Variation/Successful Run}\\
		J & 2 & 3 & 4 & 5 & 6 & 7 & 8 \\
		\hline
		1a & 0.28203 & 0.20392 & 0.081031 & 0.057056 & 0.027255 & 0.018172 & 0.016179\\
		1b & 0.31030 & 0.17807 & 0.098884 & 0.044555  & 0.028840 & 0.020392 & 0.0079454\\
		1c & 0.60090 & 0.15673 & 0.10368 & 0.054715 & 0.019420 & 0.019175 & 0.011833\\
		\hline
	\end{tabular}
	\caption{Test result for Example 1 with explicit scheme.} \label{test_result_explicit}
\end{table}

While we have not shown the proof of convergence for the Crank–Nicolson scheme, where $\theta_1 = \theta_2 = 0.5$, our numerical tests for test cases 1d, 1e, 1f, in Table \ref{test_result_crank_nicolson}, suggest that it works well in our framework.

\begin{table}[h]
	\centering
	\begin{tabular}{c|ccccccc}
		\multicolumn{8}{c}{$|Y_0-y^{(\theta_1, \theta_2), R}_0(x_0)|$}\\
		J & 2 & 3 & 4 & 5 &6 & 7 & 8 \\
		\hline
		1a & $0.0053401$ & $0.032606$ & $0.18142$ & $0.025799$ & $0.0060404$ & $0.020565$ & NA \\
		1b & 3.6788 & 0.24551 & 0.34892 & 0.069220 & 0.012861 &	0.0013653 & 0.0024095\\
		1c & $4.6822 \times 10^8$ & $3.5241 \times 10^{137}$ & $1.0773 \times 10^{44}$ & 0.051122  & 0.0050518 & 0.011735 & 0.0030526\\
		\hline
		\multicolumn{8}{c}{Total Variation/Successful Run}\\
		J & 2 & 3 & 4 & 5 & 6 & 7 & 8 \\
		\hline
		1a & $0.23450$ & $0.032606$ & $0.18142$ & $0.025799$ & $0.012630$ & $0.020565$ & NA\\
		1b & 4.5732 & 0.37590 & 0.34892 & 0.075550 & 0.014571 & 0.012470 & 0.010903\\
		1c & $4.6822 \times 10^8$ & $3.5241 \times 10^{137}$ & $1.0773 \times 10^{44}$ &	0.058288 & 0.020924 & 0.014260 & 0.0078873\\
		\hline
	\end{tabular}
	\caption{Test result for Example 1 with Crank–Nicolson scheme} \label{test_result_crank_nicolson}
\end{table}

A specific point of interest is the impact of factor $L$ introduced in Section \ref{section_refined_regression} for the samples selection.
It can be seen in Table \ref{test_result_crank_nicolson} that when the number of paths or the bundles are few, a smaller value of $L$ preserves the stability of our algorithm.
In test case $1d$, where the factor $L$ is relatively small, our algorithm rejected all tests for $J=8$.
One of the explanations is that the regression coefficients converge to the analytic projection coefficients on the basis space but the norm of these analytic coefficient is greater than $L$.
The effect of the factor $L$ actually can be seen in Table \ref{test_result_explicit} too. 
Some runs for test case 1a were rejected when $J=8$ and the result for $J=8$ is worse than either 1b or 1c. 
On the contrary, if we remove the restriction on $L$ altogether, the results are non-satisfactory when the value of $J$ is low but converge when the number of time steps and samples are high enough.
Heuristically, the regression coefficients should converge to the actual projection coefficients on the basis space, which results in a function that is bounded in a compact set.
This in turns satisfies the conditions of the proof of convergence with respect to the regression.
Although it may look like we can adjust $L$ in the same time as other algorithm parameters in order to achieve the optimal result, we should still note that $L$ is model dependent and there is no clear way to figure out the best link of $L$ with the simulation parameters. 
It remains important to use $L$ as a warning system.

Next, we shall move on to the result for the more practical and higher-dimensional Example 2.
The results for Example 2.1 in Table \ref{test_result_arithmetic} show that our method can be easily applied to a practical problem.

\begin{table}[h]
	\centering
	\centering
	\begin{tabular}{c|ccccc}
		\multicolumn{6}{c}{$|Y_0-y^{(\theta_1, \theta_2), R}_0(x_0)|$}\\
		J & 0 & 1 & 2 &  &  \\
		\hline
		2.1a & $2.0321 \times 10^{-3}$ & $2.2567 \times 10^{-3}$ & $1.9883 \times 10^{-3}$ &  & \\
		2.1b & $2.9314 \times 10^{-3}$ & $1.8934 \times 10^{-3}$ & $2.2151 \times 10^{-4}$ &  & 
	\end{tabular}
	\caption{Test result for Example 2.1 } \label{test_result_arithmetic}
\end{table}

With respect to the problem of dimensionality, we can check the results in Table \ref{test_result_geometric}.
Since the analytic solution is known to this problem, we compare our result to the actual value.
It can be seen that under our choice of bundling and regression basis, the accuracy of our method is similar across all choices of problem dimensions.
This suggested that with appropriate setting, our algorithm can easily scale up to tackle high-dimensional problems.

\begin{table}[h]
	\centering
	\begin{tabular}{c|ccccc}
		\multicolumn{6}{c}{$|Y_0-y^{(\theta_1, \theta_2), R}_0(x_0)|$}\\
		Stock dimensions & 1 & 2 & 3 & 4 & 5 \\
		\hline
		2.2a & $6.5482 \times 10^{-3}$ & $7.3015 \times 10^{-3}$ & $6.6827 \times 10^{-3}$ & $8.0384 \times 10^{-3}$ & $7.1308 \times 10^{-3}$\\
		2.2b & $5.1918 \times 10^{-3}$ & $6.9460 \times 10^{-3}$ & $6.4038 \times 10^{-3}$ & $6.9507 \times 10^{-3}$ & $7.4937 \times 10^{-3}$\\
		\hline
		Stock dimensions & 6 & 7 & 8 & 9 & 10\\
		\hline
		2.2a & $6.9885 \times 10^{-3}$ & $7.5067 \times 10^{-3}$ & $6.9271 \times 10^{-3}$ & $6.9993 \times 10^{-3}$ & $7.5682 \times 10^{-3}$\\
		2.2b & $7.2034 \times 10^{-3}$ & $7.1633 \times 10^{-3}$ & $7.0850 \times 10^{-3}$ & $7.2023 \times 10^{-3}$ & $6.7595 \times 10^{-3}$\\
		\hline
		Stock dimensions & 11 & 12 & 13 & 14 & 15\\
		\hline
		2.2a & $6.9549 \times 10^{-3}$ & $7.4005 \times 10^{-3}$ & $7.5329 \times 10^{-3}$ & $7.1437 \times 10^{-3}$ & $7.1364 \times 10^{-3}$\\
		2.2b & $8.4614 \times 10^{-3}$ & $7.1430 \times 10^{-3}$ & $7.6267 \times 10^{-3}$ & $7.8998 \times 10^{-3}$ & $7.2455 \times 10^{-3}$\\
	\end{tabular}
	\caption{Test result for Example 2.2} \label{test_result_geometric}
\end{table}

More generally, all the results from Example 2 suggest that linking the bundling criterion and the regression basis to the terminal condition can deliver an accurate algorithm.
Adapting our algorithm to a specific problem to improve the performance could be a promising direction of further research.
In fact, the choice of basis itself deserves further study. 
Even in our localised setting, regression with respect to the linear basis scheme fails to converge for Example 1.
A more sophisticated way to pick the regression basis may be important to put our algorithm into actual applications.

To sum up, we have developed a new algorithm for approximating BSDEs based on SGBM and our numerical tests showed that this new algorithm can deliver accurate estimation results.

\section*{Acknowledgements}
The authors would like to thank VORtech, BV, for their help and advice for this work and the anonymous reviewers for their valuable advice for improving this work.
This work is supported by EU Framework Programme for Research and Innovation Horizon 2020 (H2020-MSCA-ITN-2014, Project 643045, 'EID WAKEUPCALL').

\bibliographystyle{plain}
\bibliography{../../SGBM}

\end{document}